\documentclass{amsart}

\input xy
\usepackage{amssymb,amsmath,stmaryrd,mathrsfs}
\def\definetac{\newif\iftac}    
\ifx\tactrue\undefined
  \definetac
  \ifx\state\undefined\tacfalse\else\tactrue\fi
\fi
\iftac\else\usepackage{amsthm}\fi
\usepackage[all,2cell]{xy}
\UseAllTwocells
\usepackage{enumitem}
\usepackage{color}
\definecolor{darkgreen}{rgb}{0,0.45,0} 
\usepackage[pagebackref,colorlinks,citecolor=darkgreen,linkcolor=darkgreen]{hyperref}
\usepackage{mathtools}          
\usepackage{braket}             

\usepackage{url}                
\usepackage{xspace}             

\makeatletter
\let\ea\expandafter

\def\mdef#1#2{\ea\ea\ea\gdef\ea\ea\noexpand#1\ea{\ea\ensuremath\ea{#2}\xspace}}
\def\alwaysmath#1{\ea\ea\ea\global\ea\ea\ea\let\ea\ea\csname your@#1\endcsname\csname #1\endcsname
  \ea\def\csname #1\endcsname{\ensuremath{\csname your@#1\endcsname}\xspace}}

\DeclareRobustCommand\widecheck[1]{{\mathpalette\@widecheck{#1}}}
\def\@widecheck#1#2{%
    \setbox\z@\hbox{\m@th$#1#2$}%
    \setbox\tw@\hbox{\m@th$#1%
       \widehat{%
          \vrule\@width\z@\@height\ht\z@
          \vrule\@height\z@\@width\wd\z@}$}%
    \dp\tw@-\ht\z@
    \@tempdima\ht\z@ \advance\@tempdima2\ht\tw@ \divide\@tempdima\thr@@
    \setbox\tw@\hbox{%
       \raise\@tempdima\hbox{\scalebox{1}[-1]{\lower\@tempdima\box
\tw@}}}%
    {\ooalign{\box\tw@ \cr \box\z@}}}

\newcount\foreachcount

\def\foreachletter#1#2#3{\foreachcount=#1
  \ea\loop\ea\ea\ea#3\@alph\foreachcount
  \advance\foreachcount by 1
  \ifnum\foreachcount<#2\repeat}

\def\foreachLetter#1#2#3{\foreachcount=#1
  \ea\loop\ea\ea\ea#3\@Alph\foreachcount
  \advance\foreachcount by 1
  \ifnum\foreachcount<#2\repeat}

\def\definescr#1{\ea\gdef\csname s#1\endcsname{\ensuremath{\mathscr{#1}}\xspace}}
\foreachLetter{1}{27}{\definescr}
\def\definecal#1{\ea\gdef\csname c#1\endcsname{\ensuremath{\mathcal{#1}}\xspace}}
\foreachLetter{1}{27}{\definecal}
\def\definebold#1{\ea\gdef\csname b#1\endcsname{\ensuremath{\mathbf{#1}}\xspace}}
\foreachLetter{1}{27}{\definebold}
\def\definebb#1{\ea\gdef\csname l#1\endcsname{\ensuremath{\mathbb{#1}}\xspace}}
\foreachLetter{1}{27}{\definebb}
\def\definefrak#1{\ea\gdef\csname f#1\endcsname{\ensuremath{\mathfrak{#1}}\xspace}}
\foreachletter{1}{9}{\definefrak} 
\foreachletter{10}{27}{\definefrak}
\def\definebar#1{\ea\gdef\csname #1bar\endcsname{\ensuremath{\overline{#1}}\xspace}}
\foreachLetter{1}{27}{\definebar}
\foreachletter{1}{8}{\definebar} 
\foreachletter{9}{15}{\definebar} 
\foreachletter{16}{27}{\definebar}
\def\definetil#1{\ea\gdef\csname #1til\endcsname{\ensuremath{\widetilde{#1}}\xspace}}
\foreachLetter{1}{27}{\definetil}
\foreachletter{1}{27}{\definetil}
\def\definehat#1{\ea\gdef\csname #1hat\endcsname{\ensuremath{\widehat{#1}}\xspace}}
\foreachLetter{1}{27}{\definehat}
\foreachletter{1}{27}{\definehat}
\def\definechk#1{\ea\gdef\csname #1chk\endcsname{\ensuremath{\widecheck{#1}}\xspace}}
\foreachLetter{1}{27}{\definechk}
\foreachletter{1}{27}{\definechk}
\def\defineul#1{\ea\gdef\csname u#1\endcsname{\ensuremath{\underline{#1}}\xspace}}
\foreachLetter{1}{27}{\defineul}
\foreachletter{1}{27}{\defineul}

\def\autofmt@n#1\autofmt@end{\mathrm{#1}}
\def\autofmt@b#1\autofmt@end{\mathbf{#1}}
\def\autofmt@l#1#2\autofmt@end{\mathbb{#1}\mathsf{#2}}
\def\autofmt@c#1#2\autofmt@end{\mathcal{#1}\mathit{#2}}
\def\autofmt@s#1#2\autofmt@end{\mathscr{#1}\mathit{#2}}
\def\autofmt@f#1\autofmt@end{\mathfrak{#1}}
\def\autofmt@u#1\autofmt@end{\underline{\smash{\mathsf{#1}}}}
\def\autofmt@U#1\autofmt@end{\underline{\underline{\smash{\mathsf{#1}}}}}
\def\autofmt@h#1\autofmt@end{\widehat{#1}}
\def\autofmt@r#1\autofmt@end{\overline{#1}}
\def\autofmt@t#1\autofmt@end{\widetilde{#1}}
\def\autofmt@k#1\autofmt@end{\check{#1}}

\def\auto@drop#1{}
\def\autodef#1{\ea\ea\ea\@autodef\ea\ea\ea#1\ea\auto@drop\string#1\autodef@end}
\def\@autodef#1#2#3\autodef@end{%
  \ea\def\ea#1\ea{\ea\ensuremath\ea{\csname autofmt@#2\endcsname#3\autofmt@end}\xspace}}
\def\autodefs@end{blarg!}
\def\autodefs#1{\@autodefs#1\autodefs@end}
\def\@autodefs#1{\ifx#1\autodefs@end%
  \def\autodefs@next{}%
  \else%
  \def\autodefs@next{\autodef#1\@autodefs}%
  \fi\autodefs@next}

\DeclareSymbolFont{bbold}{U}{bbold}{m}{n}
\DeclareSymbolFontAlphabet{\mathbbb}{bbold}

\mdef\delbar{\overline{\partial}}

\mdef\hf{\textstyle\frac12 }
\mdef\thrd{\textstyle\frac13 }
\mdef\qtr{\textstyle\frac14 }

\newcommand{\op}{^{\mathrm{op}}}

\SelectTips{cm}{}
\newdir{ >}{{}*!/-10pt/@{>}}    

\mdef\Id{\mathrm{Id}}
\mdef\id{\mathrm{id}}
\alwaysmath{ell}
\alwaysmath{infty}
\alwaysmath{odot}
\def\frc#1/#2.{\frac{#1}{#2}}   
\mdef\ten{\mathrel{\otimes}}

\mdef\sqten{\mathrel{\boxtimes}}

\DeclareRobustCommand\widecheck[1]{{\mathpalette\@widecheck{#1}}}
\def\@widecheck#1#2{%
    \setbox\z@\hbox{\m@th$#1#2$}%
    \setbox\tw@\hbox{\m@th$#1%
       \widehat{%
          \vrule\@width\z@\@height\ht\z@
          \vrule\@height\z@\@width\wd\z@}$}%
    \dp\tw@-\ht\z@
    \@tempdima\ht\z@ \advance\@tempdima2\ht\tw@ \divide\@tempdima\thr@@
    \setbox\tw@\hbox{%
       \raise\@tempdima\hbox{\scalebox{1}[-1]{\lower\@tempdima\box
\tw@}}}%
    {\ooalign{\box\tw@ \cr \box\z@}}}

\DeclareMathOperator\Hom{Hom}

\mdef\we{\overset{\sim}{\longrightarrow}}
\mdef\leftwe{\overset{\sim}{\longleftarrow}}

\let\xto\xrightarrow

\def\rightarrowtailfill@{\arrowfill@{\Yright\joinrel\relbar}\relbar\rightarrow}
\newcommand\xrightarrowtail[2][]{\ext@arrow 0055{\rightarrowtailfill@}{#1}{#2}}

\def\twoheadrightarrowfill@{\arrowfill@{\relbar\joinrel\relbar}\relbar\twoheadrightarrow}
\newcommand\xtwoheadrightarrow[2][]{\ext@arrow 0055{\twoheadrightarrowfill@}{#1}{#2}}

\def\slashedarrowfill@#1#2#3#4#5{%
  $\m@th\thickmuskip0mu\medmuskip\thickmuskip\thinmuskip\thickmuskip
   \relax#5#1\mkern-7mu%
   \cleaders\hbox{$#5\mkern-2mu#2\mkern-2mu$}\hfill
   \mathclap{#3}\mathclap{#2}%
   \cleaders\hbox{$#5\mkern-2mu#2\mkern-2mu$}\hfill
   \mkern-7mu#4$%
}
\def\rightslashedarrowfill@{%
  \slashedarrowfill@\relbar\relbar\mapstochar\rightarrow}
\newcommand\xslashedrightarrow[2][]{%
  \ext@arrow 0055{\rightslashedarrowfill@}{#1}{#2}}
\mdef\hto{\xslashedrightarrow{}}
\mdef\htoo{\xslashedrightarrow{\quad}}

\def\shvar#1#2{{\ensuremath{%
  \hspace{1mm}\makebox[-1mm]{$#1\langle$}\makebox[0mm]{$#1\langle$}\hspace{1mm}%
  {#2}%
  \makebox[1mm]{$#1\rangle$}\makebox[0mm]{$#1\rangle$}%
}}}
\def\sh{\shvar{}}

\long\def\my@drawfill#1#2;{%
\@skipfalse
\fill[#1,draw=none] #2;
\@skiptrue
\draw[#1,fill=none] #2;
}
\newif\if@skip
\newcommand{\skipit}[1]{\if@skip\else#1\fi}
\newcommand{\drawfill}[1][]{\my@drawfill{#1}}

\newif\ifhyperref
\@ifpackageloaded{hyperref}{\hyperreftrue}{\hyperreffalse}
\iftac
  \let\your@state\state
  \def\state#1{\gdef\currthmtype{#1}\your@state{#1}}
  \let\your@staterm\staterm
  \def\staterm#1{\gdef\currthmtype{#1}\your@staterm{#1}}
  \let\defthm\newtheorem
  \def\currthmtype{}
  \ifhyperref
    \def\autoref#1{\ref*{label@name@#1}~\ref{#1}}
  \else
    \def\autoref#1{\ref{label@name@#1}~\ref{#1}}
  \fi
  \AtBeginDocument{%
    \let\old@label\label%
    \def\label#1{%
      {\let\your@currentlabel\@currentlabel%
        \edef\@currentlabel{\currthmtype}%
        \old@label{label@name@#1}}%
      \old@label{#1}}
  }
\else
  \ifhyperref
    \def\defthm#1#2{%
      \newtheorem{#1}{#2}[section]%
      \expandafter\def\csname #1autorefname\endcsname{#2}%
      \expandafter\let\csname c@#1\endcsname\c@thm}
  \else
    \def\defthm#1#2{\newtheorem{#1}[thm]{#2}}
    \ifx\SK@label\undefined\let\SK@label\label\fi
    \let\old@label\label
    \let\your@thm\@thm
    \def\@thm#1#2#3{\gdef\currthmtype{#3}\your@thm{#1}{#2}{#3}}
    \def\currthmtype{}
    \def\label#1{{\let\your@currentlabel\@currentlabel\def\@currentlabel%
        {\currthmtype~\your@currentlabel}%
        \SK@label{#1@}}\old@label{#1}}
    \def\autoref#1{\ref{#1@}}
  \fi
\fi

\newtheorem{thm}{Theorem}[section]

\defthm{cor}{Corollary}
\defthm{prop}{Proposition}
\defthm{lem}{Lemma}
\defthm{sch}{Scholium}
\defthm{assume}{Assumption}
\defthm{claim}{Claim}
\defthm{conj}{Conjecture}
\defthm{hyp}{Hypothesis}
\defthm{fact}{Fact}
\iftac\theoremstyle{plain}\else\theoremstyle{definition}\fi
\defthm{defn}{Definition}
\defthm{notn}{Notation}
\iftac\theoremstyle{plain}\else\theoremstyle{remark}\fi
\defthm{rmk}{Remark}
\defthm{eg}{Example}
\defthm{egs}{Examples}
\defthm{ex}{Exercise}
\defthm{ceg}{Counterexample}

\def\thmqedhere{\expandafter\csname\csname @currenvir\endcsname @qed\endcsname}

\setitemize[1]{leftmargin=2em}
\setenumerate[1]{leftmargin=*}

\iftac
  \let\c@equation\c@subsection
\else
  \let\c@equation\c@thm
\fi
\numberwithin{equation}{section}

\@ifpackageloaded{mathtools}{\mathtoolsset{showonlyrefs,showmanualtags}}{}

\alwaysmath{alpha}
\alwaysmath{beta}
\alwaysmath{gamma}
\alwaysmath{Gamma}
\alwaysmath{delta}
\alwaysmath{Delta}
\alwaysmath{epsilon}
\mdef\ep{\varepsilon}
\alwaysmath{zeta}
\alwaysmath{eta}
\alwaysmath{theta}
\alwaysmath{Theta}
\alwaysmath{iota}
\alwaysmath{kappa}
\alwaysmath{lambda}
\alwaysmath{Lambda}
\alwaysmath{mu}
\alwaysmath{nu}
\alwaysmath{xi}
\alwaysmath{pi}
\alwaysmath{rho}
\alwaysmath{sigma}
\alwaysmath{Sigma}
\alwaysmath{tau}
\alwaysmath{upsilon}
\alwaysmath{Upsilon}
\alwaysmath{phi}
\alwaysmath{Pi}
\alwaysmath{Phi}
\mdef\ph{\varphi}
\alwaysmath{chi}
\alwaysmath{psi}
\alwaysmath{Psi}
\alwaysmath{omega}
\alwaysmath{Omega}

\makeatother

\xyoption{all}
\usepackage{notation}

\usepackage{amsmath,amssymb,amsthm}
\usepackage{graphics}
\usepackage{enumerate}
\usepackage{mathtools}

\newcommand{\NH}{N\!H}
\newcommand{\WH}{W\!H}

\newcommand{\WK}{W\!K}
\CompileMatrices

\theoremstyle{theorem}
\newtheorem*{thmnonumber}{Theorem}
\newtheorem*{pthmA}{Preliminary Version of Theorem A}
\newtheorem*{pthmB}{Preliminary Version of Theorem B}
\newtheorem*{thmA}{Theorem A}
\newtheorem*{thmB}{Theorem B}

\let\xto\xrightarrow

\newcommand{\ecomcat}[2]{\Pi_0({#1},{#2})}
\newcommand{\ecomcatf}[2]{\Pi_0^{f}({#1},{#2})}
\newcommand{\ecomsp}[2]{\overline{{#1}|{#2}}}
\newcommand{\ecommap}[1]{\overline{{#1}}}
\newcommand{\efuncat}[2]{\Pi({#1},{#2})}
\newcommand{\efuncatf}[2]{\Pi^{f}({#1},{#2})}

\newcommand{\cen}[2]{C_{#2}({#1})}
\newcommand{\con}[2]{Conj_{#2}({#1})}
\newcommand{\rcov}[1]{\overline{#1}}
\newcommand{\gTop}{{\sf GpTop}}
\newcommand{\Iso}[1]{{\sf Iso}({#1})}
\newcommand{\Isog}[2]{{\sf Iso}({#1})({#2})}

\DeclareMathOperator{\Fix}{Fix}

\theoremstyle{plain}
\makeatletter
\newtheorem*{rep@theorem}{\rep@title}
\newcommand{\newreptheorem}[2]{%
\newenvironment{rep#1}[1]{%
 \def\rep@title{#2 \ref{##1}}%
 \begin{rep@theorem}}%
 {\end{rep@theorem}}}
\makeatother

\newreptheorem{thm}{Theorem}
\newreptheorem{lem}{Lemma}

\begin{document}

\title{Equivariant Fixed Point Theory}
\author{Kate Ponto}
\date{\today\\ The author was partially supported by NSF grants DMS-0703574 and DMS-1207670.}

\maketitle
\begin{abstract}
  We reexamine equivariant generalizations of the Lefschetz number and Reidemeister trace 
  using categorical traces.  This gives simple, conceptual descriptions of the invariants as well as direct comparisons to 
  previously defined generalizations.  These comparisons are illuminating applications of the additivity and multiplicativity of 
  the categorical trace.
\end{abstract}
\section*{Introduction}

There are two natural ways to approach generalizations of the Lefschetz fixed point theorem and its converse.  
One follows the classical description of the Nielsen number \cite{FadellWong, Wong, Wong2}.  The alternative approach \cite{KW2}
starts from the more homotopical description of \cite{DP, KW, LMS}.  
In this paper we will compare  the equivariant invariants arising from these different starting points using formal tools that makes their connections transparent.

The homotopical starting point for equivariant fixed point theory is a pair of equivariant stable homotopy classes.  
If $G$ is a finite group, $X$ is a compact $G$-ENR or closed smooth $G$-manifold and $\{\, ,\, \}_G$ denotes equivariant stable homotopy classes of maps,
 we associate an {\bf equivariant Lefschetz number} 
  $L_G(f)\in \{S^0,S^0\}_{G}$   to every   equivariant  map $f\colon X\rightarrow X$.  The equivariant Lefschetz number of the identity map is 
the {\bf equivariant Euler characteristic} $\chi_G(X)$.

  If \[\Lambda^fX\coloneqq \{\gamma\in X^I|f(\gamma(0))=\gamma(1)\}\]
is the $f$-twisted loops in $X$ we also have  an {\bf equivariant Reidemeister trace} $R_G(f)\in \{S^0,(\Lambda^fX)_+\}_{G}$.  For all of these invariants the same notation without a subscript indicates
the corresponding nonequivariant object.
\begin{thmnonumber}\cite{DP, KW2} \label{kwconverse} 
   If  $f$ has no fixed points $L_G(f)$ and $R_G(f)$ are trivial.  

  Additionally suppose $X$ is a closed smooth $G$-manifold and  for all isotropy 
  subgroups $K\subset H$ of $X$ $\mathrm{dim}(X^H)\geq 3$ and 
  $\mathrm{dim}(X^H)\leq \mathrm{dim}(X^K)-2$.   Then $f$ is equivariantly homotopic to a
  map with no fixed points if and only if $R_G(f)$ is trivial.
\end{thmnonumber}

Starting with the Nielsen number  the more natural approach is to look at classical invariants on 
isotropy subspaces $X^H\coloneqq \{x\in X\mid xh=x\text{ for all }h\in H\}$.  An equivariant map induces
a map $f^H\colon X^H\to X^H$ that is equivariant with respect to the action of the Weyl group $\WH\coloneqq \NH/H$.  After 
forgetting the $\WH$ action (and ignoring the parts of $X^H$ where $\WH$ does not act freely) we have a nonequivariant 
Lefschetz number $L(f_H)$.  Alternatively, if the $\WH$ action is cellular, we can consider the Hattori-Stallings trace $\tr_{/\WH}$ (\autoref{dualdef}, \cite{stallings}) of the map
induced on the rational cellular chain complex as a map of modules over $\WH$.  We have analogous options for the Reidemeister trace.  
In this paper we compare the resulting invariants.  In the following two results
 $G$ is a finite group, $X$ is a closed smooth $G$-manifold, 
and $f\colon X\to X$ is an equivariant endomorphism. 

\begin{pthmA}\label{thm:main}
If $\con{G}{}$ is a set of representatives for the conjugacy classes of subgroups of $G$  \[ L_G(f)= \sum_{H\in \con{G}{}} \chi_G(G/H)\frac{L(f_H)}{\chi(\WH)}.\]
\end{pthmA}

\begin{pthmB}
  There are maps $\mu$ and $\xi$ so that 
  \begin{align}
    L_G(f)=\sum_{H\in \con{G}{}}\chi_G(G/H)\sum_i(-1)^i\mu(\tr_{/\WH}(C_i(f_H;\mathbb{Q})))
  \end{align} and 
   $R_G(f)$ can be determined from \[\left\{\sum_i(-1)^i\xi(\tr_{/(\pi_1(X^H)\rtimes WH)}C_i(\tilde{f}_H;\mathbb{Q}))\right\}_{H\in \con{G}{}}.\]
\end{pthmB}

Later in the paper (pages \pageref{afree} and \pageref{facgenrei}) we will be more precise about the maps and give more explicit descriptions of the traces.  We will also give more 
symmetric formulations of the theorems. 

Various parts of this theorem can be found in the literature
\cite{FadellWong,Lait, LR, ulrich, Weber, Weber2, Wong, Wong2}, but we provide a very different proof where much of the hard work is 
outsourced to formal results for monoidal categories and bicategories \cite{DP, thesis, higher, PS_mult,PS6}.  This gives an especially transparent 
approach that applies to both the Lefschetz number and Reidemeister trace.  It is a good example of the advantages of the 
formal approach.

After briefly recalling the categorical preliminaries, we give proofs of these results in the case that the action in free.  Here we make significant 
use of ideas from the proof of the multiplicativity of traces \cite{PS_mult}.  In the last three sections we extend to the general case 
 and build on the linearity of traces \cite{PS6}.  The underlying ideas are the same in the second case but the bookkeeping is more complicated.

\begin{rmk}
The different invariants considered here have different natural generalities.  Some make 
sense for compact Lie groups, others infinite discrete groups.  In the interest of consistent 
hypotheses, we will always assume that $G$ is a finite group but it is useful to remember
that this is more restrictive than necessary in some cases. 
\end{rmk}

\subsection*{Acknowledgments} I would like to thank Mohammed Abouzaid, Frank Connolly, Peter May, 
Gun Sunyeekhan, and Bruce Williams for many helpful conversations.

\section{Duality and trace in symmetric monoidal categories}\label{symmetricmonoidal}

The trace in symmetric monoidal categories is a generalization of
the trace in linear algebra that retains many of the important properties.  In particular, it
satisfies a generalization of invariance of basis and is functorial.  The generalized
trace is a trace for endomorphisms of modules over a commutative ring, 
endomorphisms of chain complexes of modules over a commutative ring, and endomorphisms
of closed smooth manifolds or compact ENRs.  This section is a summary of \cite{DP,LMS,PS_trace}.

Let $\sV$ be a symmetric monoidal category with monoidal product $\otimes$, unit $S$, and symmetry
isomorphism $\gamma$.  

\begin{defn} An object $A$ in $\sV$ is {\bf dualizable} with {\bf dual} $B$ if there are morphisms 
  \[\eta\colon S\rightarrow A\otimes B \text{ and }\epsilon\colon B\otimes A\rightarrow S\] such that the 
  composites
   \[\xymatrix@R=3pt{A\cong S\otimes A\ar[r]^-{\eta\otimes\id}&A\otimes B\otimes A
     \ar[r]^-{\id \otimes \epsilon}&A\otimes S\cong A\\
    B\cong B\otimes S\ar[r]^-{\id\otimes \eta}&B\otimes A\otimes B
    \ar[r]^-{\epsilon \otimes \id}&S\otimes B\cong B}\] 
  are identity maps.
\end{defn}

We say a space is {\bf dualizable} if its suspension spectrum is dualizable in the stable homotopy category.  A space with a $G$ action 
is dualizable if its equivariant suspension spectrum is dualizable in the equivariant stable homotopy category.

\begin{prop}\cite[III.4.1, III.5.1]{LMS} \label{topduals}
 If $X$ is a compact $G$-ENR or closed smooth $G$-manifold 
  then $X_+\coloneqq X\amalg \ast$ is  dualizable.
\end{prop}

Surprisingly, an explicit description of the dual will not be important in this paper.

\begin{defn}
  If $A$ is dualizable with dual $B$ and $f\colon A\rightarrow A$ is an endomorphism in $\sV$, the {\bf trace} of $f$, $\tr(f)$, is the composite 
  \[\xymatrix{ S\ar[r]^-\eta&A\otimes B\ar[r]^{f\otimes \id}&A\otimes 
  B\ar[r]^\gamma&B\otimes A\ar[r]^-\epsilon&S}.\]
\end{defn}

The trace of a chain map is the alternating sum of the levelwise traces. If  $f$ is an endomorphism 
of a topological space and $H_*(-\colon \mathbb{Q})$ is the rational homology functor, the trace of $H_*(f\colon\mathbb{Q})$ is the {\bf Lefschetz number} of $f$.  
The trace of an endomorphism of  a $G$-space  in the equivariant stable homotopy 
category is the {\bf equivariant fixed point index} \cite{DP}.   

\begin{rmk}\label{namingconv}  In this paper 
  we will generally not distinguish between Lefschetz numbers (computed on the chain complex) and fixed point indices (computed on the level of spaces)
  since there are classical  identification theorems
  that show they agree in the cases of interest \cite{brownbook,DP}. 
 These identifications can be made in a way that is compatible with the approach here \cite{thesis}, further reducing the need to make these distinctions.
\end{rmk}

In the stable homotopy category and the equivariant stable homotopy category, as well as many other categories,  the trace is additive on cofiber sequences.
\begin{thm}\cite{additivity, mayadd}\label{smcadd} In a diagram of cofiber sequences 
  \[\xymatrix{A\ar[r]^-i\ar[d]^{f_A}&X\ar[r]\ar[d]^{f}&C\ar@{.>}[d]^{f_C}\\
  A\ar[r]&X\ar[r]&C
  }\] 
   where $A$ and $X$ are dualizable and the left square commutes,  $C$ is also dualizable and there is a map $f_C$ so that the remaining square commutes  
  and $\tr(f_A)+\tr(f_C)=\tr(f)$.
\end{thm}

Spaces with a group action have natural decompositions of this form.  If we let  $(H)$ denote the conjugacy class of the subgroup $H$ in 
  $G$  and, for a $G$-space $X$ and $x\in X$,  $G_x\coloneqq \{g\in G\mid xg=x\}$ then we use the notation 
 \begin{align}
  X_{(H)}\coloneqq\{x\in X\mid (G_x)=(H)\}&\hspace{.5cm} X^{(H)}\coloneqq\{x\in X\mid \exists g\in G\text{ where } gHg^{-1}\subset G_x\} \\
   X^{>(H)}&\coloneqq X^{(H)}\setminus X_{(H)}.\end{align}  
Each of the inclusion maps $X^{>(H)}\to X^{(H)}$ is a cofibration  \cite[II.1.9, II.6.7]{ulrich}.

\begin{thm}\cite[III.5.4]{ulrich}\label{ulrichadd}
  If $X$ is a closed smooth $G$-manifold or compact $G$-ENR then 
  \[L_G(f)=\sum_{H\in \con{G}{}} L_G(f_{{(H)}})\]
  where $f_{(H)}$ is the induced
  endomorphism of  $X^{(H)}/X^{>(H)}$.
\end{thm}

\begin{proof}
  Containment defines a partial order on the set of conjugacy classes of subgroups of $G$.  Extend this to a total order
  \[(e)= (H_1)< (H_2)< (H_3)< ...<  (H_n )= G.\]
  By \cite[II.6.7]{ulrich}, $X^{(H_i)}$ and $X^{>(H_i)}$ are compact $G$-ENRs and so they are dualizable.
  Then the map of cofiber sequences 
  \[\xymatrix{
    X^{>(H_i)}\ar[r]\ar[d]^{f^{>(H_i)}}& X^{(H_i)}\ar[r]\ar[d]^{f^{(H_i)}}& C_i\ar[d]^{f_i}\\
    X^{>(H_i)}\ar[r]& X^{(H_i)}\ar[r]& C_i
  }\]
  and  \autoref{smcadd} imply $L_G(f^{(H_i)})=L_G(f^{>(H_i)})+L_G(f_i)$.  In this case we 
can take $C_i=X^{(H_i)}/X_{(H_i)}$ and $f_i=f_{(H_i)}$.

Then 
$L_G(f^{>(H_k)})$ can be written as a sum of $L_G(f_{(H_i)})$ where $i\geq k$.\end{proof}

In papers such as \cite{Wilczynski}, these types of decompositions play an essential role, but they are expressed
 in terms of \emph{taut} maps.  A map $f\colon X\rightarrow Y$ is {\bf taut} if for all 
isotropy subgroups $H$ of $X$ there is a neighborhood $V$ of $X^{>H}\coloneqq\{x\in X\mid H\subsetneq G_x\}$ in $X^H\coloneqq\{x\in X\mid H=G_x\}$ and an equivariant retraction
 $r_H\colon V\rightarrow X^{>H}$ such that $f^H\rvert_V=f^H\circ r_H$.
The assumption that $X$ is a compact $G$-ENR or a closed smooth $G$-manifold implies that
any equivariant endomorphism of $X$ is equivariantly homotopic to a taut map.  Since the invariants here are all defined 
up to homotopy,  taut maps will not play an explicit role in this paper.

\section{Duality and trace in bicategories with shadows}\label{dualbicat}

To define the Reidemeister trace from this perspective and to capture the comparison results in Theorems A and B we need to extend the trace in 
a symmetric monoidal category to a bicategory.  
This section is a brief summary of the relevant parts of \cite{MS,thesis,shadows}.

\begin{defn} A {\bf bicategory} $\sB$ consists of
  \begin{itemize}\item A collection $\ob\sB$. 
     \item Categories $\sB(A,B)$ for each $A,B\in \mathrm{ob}\sB$. 
    \item Functors \[\odot \colon \sB(A,B)\times \sB(B,C)\rightarrow\sB(A,C)\]
      \[U_A\colon \ast \rightarrow \sB(A,A)\] for $A$, $B$ and $C$ in $\mathrm{ob}\sB$.
  \end{itemize}
  Here $\ast$ denotes the category with one object and one morphism.
  The functors $\odot$ are  required to satisfy unit and
  associativity axioms  up to  natural isomorphisms in $\sB(A,B)$.
\end{defn}

The elements of $\ob\sB$ are called \emph{0-cells}.  The objects of $\sB(A,B)$
are called \emph{1-cells}.  The morphisms of $\sB(A,B)$ are called \emph{2-cells}.

\begin{eg}\label{bicateg}  
  \begin{itemize}
     \item The 0-cells in the bicategory $\Mod$  are rings and 
    the category $\Mod(R,S)$ for rings $R$ and $S$ is the category of $R$-$S$-bimodules and their homomorphisms.   
    The composition is given by tensor product and a ring regarded as a module over itself is the unit.  

    \item The 0-cells in the bicategory $\Ch$  are rings and 
    the category $\Ch(R,S)$ for rings $R$ and $S$ is the category of chain complexes of $R$-$S$-bimodules and their chain homotopy classes of maps.  
    The composition is given by tensor product and a ring regarded as a module over itself is the unit.  

    \item The 0-cells in the bicategory $\gTop$ are finite groups.   
     A 1-cell $X\colon G\to H$ is a based space with an action of $G\times H$ where the actions of $G$ and $H$ are separately free away from the base point.
    The morphisms from $X\colon G\to H$ to $Y\colon G\to H$ are stable homotopy classes of equivariant maps
    from $X$ to $Y$.
    The bicategorical composition is given by the smash product followed by the quotient by the diagonal action. 
    The unit object associated to a finite group $G$ is $G_+$ regarded as a $G$-$G$ set with a trivial action on the base point.

    \item The 0-cells in the bicategory $\Ex$ of parametrized spectra defined in \cite{MS} are 
    topological spaces.  The 1-cells are parametrized spectra and the 2-cells are fiberwise stable homotopy classes of maps.
    The bicategory composition is given by a fiberwise smash product.   For this bicategory we will follow the notation and conventions of \cite[\S 3]{PS_mult}.

     \noindent There is also a bicategory of parametrized spectra with 
    an action by a finite group $G$. 

  \end{itemize}
\end{eg}

The first two of these bicategories primarily serve as motivation.  Our interest is in bicategories arising in topological settings. 

\begin{defn}\cite[16.4.1]{MS} A 1-cell $X\in \sB(A,B)$ is {\bf right dualizable} with 
  dual $Y\in \sB(B,A)$ if there are 2-cells \[\xymatrix{\eta\colon U_A\ar[r]&X\odot Y&\epsilon\colon Y\odot X
  \ar[r] &U_B}\] such that  the composites 
  \[\xymatrix@R=5pt{Y\cong Y\odot U_A\ar[r]^-{\id \odot \eta}&
  Y\odot X\odot Y\ar[r]^-{\epsilon\odot \id}&
  U_B\odot Y\cong Y\\
  X\cong U_A\odot X\ar[r]^-{\eta\odot \id}&
  X\odot Y\odot X\ar[r]^-{\id \odot \epsilon}&
  X\odot U_B\cong X}\]
  are identity maps.
\end{defn}

The map $\eta$ is the {\bf coevaluation} and $\epsilon$ is the {\bf evaluation}.  We say $(X,Y)$ is a {\bf dual pair}.

In this paper we will use a range of topological dual pairs.  Most are closely related to the classical dual pair in \autoref{topduals}.

\begin{thm}\cite[8.6]{Adams}\label{adamsdual} If $X$ is a compact ENR or closed smooth manifold with a free right action of a finite group $G$ the space $X_+$ is dualizable as a
$\ast\times G$ space 
in the bicategory $\gTop$.
\end{thm}

The evaluation and coevaluation can be 
interpreted as maps of spaces
\begin{equation}\label{randual}S^n\to X_+\wedge_G DX\text{ and  } DX\wedge X_+\to G_+\wedge S^n\end{equation}
for a $G$ space $DX$ and an integer  $n$ \cite{Adams}.  The first map is a map of spaces and  the second map is $G$-$G$-equivariant.
To give a better idea of the objects involved we use the notation $\wedge_G$ and $\wedge$ rather than the standard $\odot$ notation.  In this case $\wedge_G$ is the smash product followed by the 
quotient by the diagonal $G$ action.

\begin{thm}\cite[3.2.3]{thesis} For a closed smooth manifold or compact ENR $X$ the universal cover $\tilde{X}$
is dualizable as a $\ast\times \pi_1(X)$ space in the bicategory $\gTop$. 
\end{thm}

We say a parametrized space $E$ over $Y\times X$
is \emph{dualizable} if the fiberwise suspension spectrum $\Sigma_{Y\times X}E$ is dualizable.  

\begin{thm}\cite[18.5.1, 18.6.1]{MS}  If $X$ is a compact $G$-ENR or closed smooth $G$-manifold 
  $S^0_X\coloneqq X\amalg X$, regarded as a parametrized space over $\ast\times X$,  is right dualizable.
\label{paradual} \end{thm}

From a map of topological spaces $f\colon X\to Y$, we define spaces $P(\id,f)\coloneqq \{(\gamma,x)\in Y^I\times X|\gamma(0)=f(x) \}$ 
and $P(f,\id)\coloneqq \{(x,\gamma)\in X\times Y^I|\gamma(1)=f(x)\}.$ The first has a map to $Y\times X$ by $(\gamma,x)\mapsto (\gamma(1),x)$ and 
the second has a similar map to $X\times Y$.  These become parametrized spaces with the addition of a disjoint section.  We let $Y_f\coloneqq P(\id,f)\amalg (Y\times X)$ and ${}_fY
\coloneqq P(f,\id)\amalg (X\times Y)$.
\begin{thm}\label{basechangedual}\cite[17.3.1]{MS}
For any map of spaces $f\colon X\to Y$ $({}_fY,Y_f)$ is a dual pair.
\end{thm}
Composition of paths and applying the map $f$ to a path defines evaluation and coevaluation maps for this dual pair. This type of dual pair will be referred to as a \emph{base change dual pair} \cite[17.3]{MS}.

Like the symmetric monoidal trace, the trace of a 2-cell  is defined using a 
composite of the coevaluation and evaluation for a dual pair. 
Unlike that case,  the source of the evaluation and target of the coevaluation are
not isomorphic.   To accommodate this, we need more structure on a bicategory 
before we can define the trace.

\begin{defn}\cite[4.4.1]{thesis}
   A {\bf shadow} for a bicategory $\sB$ is a functor \[\sh{-}\colon \coprod\sB(A,A)\rightarrow \sT\] to a 
   category $\sT$ and  unital and associative natural isomorphisms $\sh{X\odot Y}\cong \sh{Y\odot X}$ for every pair of 
   1-cells $X\in \sB(A,B)$ and $Y\in \sB(B,A)$.
\end{defn}
All of the bicategories in \autoref{bicateg}  have shadows \cite{shadows}.  The shadow in $\gTop$ is the quotient by the diagonal action of the group.  In the bicategory $\Ex$ the shadow
is given by pulling back along the diagonal map (up to homotopy)  and then quotienting by the resulting section.  In particular,  for an endomorphism $f\colon X\to X$, 
$\sh{X_f}\cong (\Lambda^fX)_+$.

\begin{defn}\cite[4.5]{thesis}\label{dualdef}
  Let $X$ be a dualizable 1-cell in $\sB$ with dual $Y$ and $f\colon Q\odot X\rightarrow X\odot P$ be a 2-cell in 
  $\sB$.  The {\bf trace} of $f$ is the composite
  \[\xymatrix{{\sh{Q}}\cong \sh{Q\odot U_A}\ar[r]^-{\id \odot \eta}&
  \sh{Q\odot X\odot Y}\ar[d]^-{f\odot \id}\\
  &\sh{X\odot P\odot Y}\ar[r]^\sim&
  \sh{P\odot Y\odot X}\ar[r]^-{\id \odot \epsilon}&\sh{P\odot U_B}\cong\sh{P}.}\] 
\end{defn}

If $M$ is a finitely generated projective right $R$-module, $M$ is right dualizable and the trace of an endomorphism of $M$ 
is the {\bf Hattori-Stallings trace}.

If $G$ acts freely on a closed smooth manifold or compact ENR $X$, the trace of an equivariant map $f\colon X\to X$ with respect to the dual pair  in \autoref{adamsdual} is a map 
\[\tr_{/G}(f)\colon S^n\to S^n\wedge \sh{G_+}. \] This is another equivariant generalization of the classical fixed point index.  We will see in 
Theorem B that it is closely related to the equivariant generalization of the index defined in the previous section.

\begin{rmk}
If we apply the rational cellular chain complex functor $C_i(-;\mathbb{Q})$ to the maps \eqref{randual} we obtain a dual pair in the bicategory 
of rings, chain complexes, and homomorphisms.  An equivariant map $f\colon X\to X$ defines a map 
of chain complexes and the trace of this map 
is \[\sum_i(-1)^i\tr_{/G}(C_i(f;\mathbb{Q}))\] where $\tr_{/G}$ is the levelwise Hattori-Stallings trace. This is the \emph{universal Lefschetz class} from \cite[1.7]{Lait}. 
By functoriality of the trace \cite{shadows} under the isomorphism $\pi_0^s(\sh{G}_+)\cong 
\Hom(\mathbb{Q},\mathbb{Q}\sh{G})$ this agrees with $\tr_{/G}(f)\colon S^n\to S^n\wedge \sh{G_+}$.   
\end{rmk}

If a space $X$ has a universal cover $\tilde{X}$, an endomorphism  $f\colon X\to X$ 
defines an endomorphism $\tilde{f}$
of $\tilde{X}$ that is $\pi_1(X)$-equivariant in the sense that for $\alpha\in \pi_1(X)$ and $\tilde{x}\in \tilde{X}$
\[\tilde{f}(x\alpha )=\tilde{f}(x)f_*(\alpha)\] (with some care with base points).    
Consistent with notation
earlier we let $(\pi_1X)_{f_*+}$ be the set $(\pi_1X)_+$ with a standard left action of $\pi_1(X)$  and a right action of $\pi_1(X)$ that is first twisted by $f_*$.  We can then interpret $\tilde{f}$ as an equivariant map $\tilde{X}_+\to \tilde{X}_+\wedge _{\pi_1(X)} (\pi_1X)_{f_*+}$.  
The {\bf Reidemeister trace} of $f$,  $R(f)$, is the bicategorical trace of $\tilde{f}$ \cite{thesis}.  It is an element of the zeroth stable homotopy group of the set
 $ \sh{\pi_1X_{f_*}}\coloneqq \pi_1X/ (\gamma f_*(\delta)\sim \delta \gamma)$.

The map   $f\colon X\to X$ also 
defines a fiberwise map 
\[S^0_X\to S^0_X\odot (X_f).\]  See \cite[2.3]{higher} and \cite{crabb}.  If $X$ is a closed smooth manifold or compact ENR
the trace of this map is an element of the zeroth stable homotopy group of $\sh{X_f}\cong (\Lambda^fX)_+$.

\begin{thm}\label{idrantrace}\cite[2.3]{higher}  
  There is a natural map $\sh{X_f}\to \sh{\pi_1(X)_{f_*}}$ and the image of the trace of the fiberwise map 
  $S^0_X\xto{f} S^0_X\odot X_f$ under this map is the Reidemeister trace of $f$.
\end{thm}

As the natural map $\sh{X_f}\to \sh{\pi_1(X)_{f_*}}$ is an isomorphism on components we will following \autoref{namingconv} and refer to the trace of $S^0_X\to S^0_X\odot X_f$ as the \emph{Reidemeister trace} of $f$.

If $f\colon X\to X$ is an equivariant endomorphism of a closed smooth $G$-manifold or compact $G$-ENR 
the {\bf equivariant Reidemeister trace} of $f$ is defined to be trace of $S^0_X\to S^0_X\odot X_f$ in the equivariant parametrized stable homotopy category.  
This class is denoted is $R_G(f)$.

Like the symmetric monoidal trace, the bicategorical trace is additive.

\begin{thm}\cite{linearity, PS6} For a diagram of cofiber sequences in the parametrized stable homotopy 
  category or its equivariant generalization
  \[\xymatrix{A\ar[r]^-i\ar[d]^{f_A}&X\ar[r]\ar[d]^{f}&C\ar@{.>}[d]^{f_C}\\
  A\odot P\ar[r]&X\odot P\ar[r]&C\odot P
  }\] 
   where $A$ and $X$ are dualizable and the left square commutes,  $C$ also dualizable and there is a map $f_C$ so that the right square commutes and 
  $\tr(f_A)+\tr(f_C)=\tr(f)$.
\end{thm}

The approach used in \autoref{ulrichadd} then extends to the Reidemeister trace.  
\begin{thm}\cite[6.3]{PS6}\label{radd1} 
  If $X$ is a closed smooth $G$-manifold or compact $G$-ENR and $f\colon X\to X$ is an equivariant endomorphism then 
  \[R_G(f)=\sum_{H\in\con{G}{}}i_{(H)} R_G(f_{(H)}).\]
\end{thm}
Here $R_G(f_{(H)})$ is the {\bf equivariant relative Reidemeister trace} of $X^{(H)}$ relative to the subspace $X^{>(H)}$.  See \cite{relative} and \S\ref{hoinvgen}.  It is a refinement of
the equivariant Reidemeister trace of  $X^{(H)}/X^{>(H)}$ and takes values in $\{S^0,\Lambda^{f^{(H)}} X^{(H)}\}_G$.  
The map $i_{(H)}$ is the inclusion $\Lambda^{f^{(H)}} X^{(H)}\to \Lambda^{f} X$.

We will also use the compatibility of the trace with composites of dual pairs.  This was an essential piece of the proofs of multiplicativity \cite{PS_mult}
and additivity \cite{PS6}. 

\begin{thm} {\cite[16.5.1]{MS}}{\cite[5.4]{PS_mult}}\label{thm:compose-duality}
  If $M\in \sB( A,B)$ and $N\in \sB( B, C)$ are right dualizable, then so is $M\odot N \in\sB( A, C)$. 

  Let $Q\in \sB( A, A)$, $P\in \sB( B, B)$, and $R\in\sB( C, C)$ be 1-cells, and let $f\colon Q\odot M\to M\odot P$ and $g\colon P\odot N\to N\odot R$ be 2-cells.
  Then the following triangle commutes.
  \[ \xymatrix@C=32pt{ \sh{Q} \ar[rr]^{\tr((\id_M \odot g)\circ (f\odot\id_N))} \ar[dr]_{\tr(f)} & &
   \sh{R} \\
    & \sh{P} \ar[ur]_{\tr(g)}
    }\]
\end{thm}

\part{Free actions}

An equivariant map $f\colon X\to X$ induces
a map $\bar{f}\colon X/G\to X/G$ so the diagram below commutes.
\[\xymatrix{X\ar[r]^f\ar[d]&X\ar[d]\\X/G\ar[r]^{\bar{f}}&X/G}\]  
If the action of $G$ on $X$ is free each $\bar{f}$-twisted loop $\gamma$ in $X/G$ and lift $\tilde{\gamma}$ of $\gamma$ to $X$ define a group element $g\in G$ by  $f(\tilde{\gamma}(0))=\tilde{\gamma}(1)g$.  Up to conjugacy, this group element depends only on $\gamma$. 
If $\sh{G}$ is the  set of conjugacy classes of elements of $G$ we  define a map  
 \begin{equation}\label{groupelt}\Theta\colon \Lambda^{\bar{f}}(X/G)\to \sh{G}  \end{equation} 
by
 $\Theta(\gamma)=g$ if there is a lift $\tilde{\gamma}$ of $\gamma$ such that $f(\tilde{\gamma}(0))=\tilde{\gamma}(1)g$. 
Let $\Fix(\overline{f})(e)\coloneqq \Theta^{-1}(e)$.  
We will also let $\Theta$ denote the corresponding map $\sh{\pi_1(X/G)_{\bar{f}_*}}\to \sh{G}$

\begin{thmA}[Free case]\label{afree} If action of $G$ on $X$ is free 
then $\chi(G) L_G(f)=  \chi_G(G)L(f)$
and there are integers $a_\gamma$ so that \[R_G(f)=\sum_{\gamma\in \Fix(\bar{f})(e)} a_\gamma (i_{\gamma}\circ \tr^\triangle_G(F_\gamma))\text{ and }R(f)=\sum_{\gamma\in \Fix(\bar{f})(e)} a_\gamma (i_{\gamma}\circ \tr^\triangle(F_{\gamma})).\]
\end{thmA}

In this statement $\tr^\triangle$ is the {\bf transfer} \cite[3.7]{LMS}.  If $X$ is a dualizable object in a symmetric monoidal category 
the transfer of $X$ with respect to a map $\triangle \colon X\to X\wedge X$ is the composite
\[S\xto{\eta} X\wedge DX\xto{\triangle \wedge 1} X\wedge X\wedge DX\xto{1\wedge \gamma} X\wedge DX\wedge X\xto{1\wedge \epsilon} X\wedge S\cong X.\]  For topological spaces
with a disjoint base point we will use the map induced by the diagonal.
The transfer leads to some unavoidable asymmetry in this statement since  we cannot multiply by transfers as easily as Euler characteristics.
Here $F_\gamma$ is the fiber over $\gamma(1)$ and  the map $i_\gamma\colon F_\gamma\to \Lambda^fX$ is the inclusion of the fiber as constant paths.

For Theorem B we first need to give a more precise description of the trace we will 
compare to the Reidemeister trace.  Let $\tilde{X}$ be the universal cover of $X$.  Then $\tilde{X}$ is  a cover of $X/G$ and the action of 
$\pi_1(X/G)$ on $\tilde{X}$ encodes both the $G$ action and the action of $\pi_1(X)$.  As in the classical case, $\tilde{X}_+$ is dualizable as 
a $\pi_1(X/G)$ space.
An equivariant map $f\colon X\to X$ induces a map $\tilde{f}\colon \tilde{X}_+\to \tilde{X}_+\wedge (\pi_1(X/G)_{\bar{f}_*})_+$ and the trace of $\tilde{f}$ is a map 
\[\tr_{/\pi_1(X/G)}(\tilde{f})\colon S^n\to S^n\wedge \sh{\pi_1(X/G)_{\bar{f}_*}}_+\]

Both $\tr_{/G}(f)$ and $\tr_{/\pi_1(X/G)}(\tilde{f})$ are carrying more information than the corresponding classical invariants.  To be able to 
compare invariants we need to be able to separate out this extra information.  
For each $g\in G$, there is a map $\mu_g\colon \sh{G}_+\to S^0$ that takes all conjugacy classes of $G$ except the class that contains  $g$ to the basepoint.
There is also a map  $\xi_g\colon  \sh{\pi_1(X/G)_{\bar{f}_*}}_+\to \Theta^{-1}(g)_+$ that is the identity on $\Theta^{-1}(g)$ and 
takes all other elements to the basepoint.  Let $\zeta_g \colon \sh{\pi_1(X)_{fg_*}}_+\to \Theta^{-1}(g)_+$
be the composite $\sh{\pi_1(X)_{fg_*}}_+\to \sh{\pi_1(X/G)_{\bar{f}_*}}_+\xto{\xi_g} \Theta^{-1}(g)_+$ where the first map is induced by the 
quotient map $X\to X/G$. 

We abuse notation and let $\xi_g$  denote the composite $\left(\Lambda^{\bar{f}}X/G\right)_+\to \sh{\pi_1(X/G)_{\bar{f}_*}}_+\xto{\xi_g} \Theta^{-1}(g)_+$ and $\zeta_g$ denote the composite $\left(\Lambda^{f\cdot g}X\right)_+ \to \sh{\pi_1(X)_{fg_*}}_+\xto{\zeta_g}  \Theta^{-1}(g)_+$. 

\begin{thmB}[Free case]\label{bfree} Let $\cen{g}{G}$ be the centralizer of $g$ in $G$.
For each $g\in G$ 
  \begin{align}
    L(f\cdot g)&= \lvert \cen{g}{G}\rvert\mu_g(\tr_{/G}(f))\\
   \zeta_g R(f\cdot g)&=\lvert \cen{g}{G}\rvert   \xi_g\left(\tr_{/\pi_1(X/G)}\left(\tilde{f}\right)\right)
  \end{align} 
\end{thmB}

Beyond the simplification to the free case, this statement differs from the statement in the introduction in two ways.  
We can use functoriality of the trace \cite{thesis} to recover the algebraic descriptions of the traces.  (This is the description of the equivariant Reidemeister trace in \cite{Weber}.)
We can  use Theorem A to recover the original comparison of $\tr_{/G}$ with $L_G$.    

The essential underlying observation in this part is that the quotient map $X\to X/G$ is a covering map.  
This allows us to use the multiplicativity of the Lefschetz number and Reidemeister trace \cite{PS_mult}
to express the invariants for $X$ in terms of the invariants for $X/G$ and the fiber.
The fibers of $X\to X/G$ are finite, discrete, and isomorphic to $G$  and the endomorphisms of these fibers
induced by an equivariant map of $X$ are determined by their value on a single point.  This allows us to 
easily compute the \emph{fiberwise Lefschetz number} and \emph{fiberwise Reidemeister trace} of maps $f\colon X\to X$ 
regarded as maps over $X/G$. 

We prove Theorem A in \S\ref{hoinv}.  We prove Theorem B for the Lefschetz number in \S\ref{leffree}
and for the Reidemeister trace in \S\ref{reifree}.

\section{Homotopical invariants} \label{hoinv}

We follow the notation of \cite{PS_mult} and use $(\widehat{L_{X/G}})(f)$ and $(\widehat{R_{X/G}})(f)$ to denote the fiberwise Lefschetz number and Reidemeister traces.  These 
are stable maps 
\[(\widehat{L_{X/G}})(f)\colon \Lambda^{\bar{f}}(X/G)\to S^0\text{ and } (\widehat{R_{X/G}})(f)\colon \Lambda^{\bar{f}}(X/G)\to \Lambda^f(X).\]
  Given a path $\gamma$ in $X/G$ from $x$ to $\bar{f}(x)$ we define an endomorphism of the fiber over $\bar{f}(x)$ 
  \begin{equation}\label{fibendo}
  F_{\bar{f}(x)}\to F_x\xto{f}F_{\bar{f}(x)}
  \end{equation}
  where the first map is induced by the path $\gamma$.  Then $(\widehat{L_{X/G}})(f)(\gamma)$
  is the Lefschetz number of \eqref{fibendo} and $(\widehat{R_{X/G}})(f)(\gamma)$ is the 
 Reidemeister trace of \eqref{fibendo} composed with the inclusion of the \eqref{fibendo}
  twisted loops in $F_{\bar{f}(x)}$ into the $f$ twisted loops in $X$. 
\cite{PS_mult} considers only the nonequivariant case, but the same approach immediately generalizes 
to  equivariant invariants.  We denote these invariants by   $(\widehat{L_{X/G}})_G$
 and $(\widehat{R_{X/G}})_G$.

\begin{prop}\label{fibcompute}
  For an equivariant map $f\colon X\to X$ and a $\bar{f}$-twisted loop $\gamma $ in $X/G$ 
  \[(\widehat{L_{X/G}})(f)(\gamma)=\left\lbrace\begin{array}{ll}\chi(G)&\text{if }\Theta(\gamma)=e\\
  0&\text{otherwise}\end{array}\right.
  \]\[
  (\widehat{L_{X/G}})_G(f)(\gamma)=\left\lbrace\begin{array}{ll}\chi_G(G)&\text{if }\Theta(\gamma)=e\\
  0&\text{otherwise}\end{array}\right.\]
  \[(\widehat{R_{X/G}})(f)(\gamma)=\left\lbrace\begin{array}{ll}i_\gamma \tr^\triangle(F_{\gamma})&\text{if }\Theta(\gamma)=e\\
  0&\text{otherwise}\end{array}\right.
  \]\[
  (\widehat{R_{X/G}})_G(f)(\gamma)=\left\lbrace\begin{array}{ll}i_\gamma\tr^\triangle_G(F_{\gamma})&\text{if }\Theta(\gamma)=e\\
  0&\text{otherwise}\end{array}\right.\]
\end{prop}

\begin{proof}
  For a lift $\tilde{\gamma}$ of $\gamma$, the image of $\tilde\gamma(1)$ under the first map in \eqref{fibendo} is $\tilde{\gamma}(0)$ and so the image
  under \eqref{fibendo} is $f(\tilde\gamma(0))=\tilde{\gamma}(1)g$ for some $g\in G$.
  Since the endomorphism in \eqref{fibendo} is equivariant and $F_{\bar{f}(x)}$ is $G$-isomorphic to the $G$-set $G$,  
  the endomorphism is determined by the image of one point and so \eqref{fibendo} is multiplication by $g$.  In particular, \eqref{fibendo} 
  is the identity map if $g$ is the identity element of $G$ and has no fixed points if $g$ is not the identity element of $G$.

  The Euler characteristic of $F_{\bar{f}(x)}$ is the same as the Euler characteristic of $G$ so
   $(\widehat{L_{X/G}})(f)(\gamma)=\chi(G)$ and $(\widehat{L_{X/G}})_G(f)(\gamma)=\chi_G(G)$ if $\Theta(\gamma)=e$ \cite[1.9, 6.6]{PS_mult}.  Both invariants are 
  zero if $\Theta(\gamma)\neq e$.

  The Reidemeister trace of the identity map of a discrete space is the transfer and the Reidemeister trace of a map with no 
  fixed points is zero.  Then 
$(\widehat{R_{X/G}})(f)(\gamma) =  i_\gamma \tr^\triangle(F_{\gamma})$ and $  (\widehat{R_{X/G}})_G(f)(\gamma)=  i_\gamma\tr^\triangle_G(F_{\gamma})$ if $\Theta(\gamma)=e$
  \cite[1.13, 7.6]{PS_mult}.  
  Both invariants are zero otherwise. 
\end{proof}

\begin{proof}[Proof of Theorem A (Free case)]
  In \cite[1.17]{PS_mult}, we showed the triangles 
  \begin{equation}\label{fib_reidemeister}
  \xymatrix@C=40pt{
  S^0\ar[r]^-{R(\bar{f})}\ar[dr]_{L(f)}& \sh{(X/G)_{\bar{f}}}\ar[d]^{\widehat{L_{X/G}}(f)}&
  S^0\ar[r]^-{R(\bar{f})}\ar[dr]_{R(f)}& \sh{(X/G)_{\bar{f}}}\ar[d]^{\widehat{R_{X/G}}(f)}
  \\
  &S^0,&&\sh{X_f}
  }\end{equation}
  commute.  The same approach also  shows the corresponding 
  equivariant generalizations commute.  (We use the trivial $G$ action on $X/G$
  and so if we regard $R(f)$ as an element of $\{S^0, \Lambda^{\bar{f}} X/G\}_G$ using the trivial action we have $R_G(\bar{f})=R(\bar{f})$.)

  By \autoref{fibcompute}  the values of $\chi(G)\cdot  (\widehat{L_{X/G}})_G(f)$ and $\chi_G(G)\cdot  (\widehat{L_{X/G}})(f)$ on any twisted loop
  in $X/G$ are the same  and so 
  \begin{align*}
  \chi(G)\cdot L_G(f)=\chi(G)\cdot  (\widehat{L_{X/G}})_G(f)\circ R(\bar{f})
  &=
  \chi_G(G)\cdot  (\widehat{L_{X/G}})(f)\circ R(\bar{f})\\
  &=\chi_G(G)\cdot L(f).
  \end{align*}
  Using the isomorphism $\{S^0, \sh{(X/G)_{\bar{f}}}_+\}\cong \mathbb{Z}\pi_0( \sh{(X/G)_{\bar{f}}}_+)$
we can write $R(\bar{f})$ as $\displaystyle\sum a_\gamma \gamma $ where $a_\gamma$ are integers and  $\gamma \in \pi_0( \sh{(X/G)_{\bar{f}}}_+)$. 
Then the Reidemeister trace result  follows from the second diagram in \eqref{fib_reidemeister} and \autoref{fibcompute}.
\end{proof}

For later results we will need to know a little more about the fiberwise Lefschetz number and Reidemeister trace.
\begin{prop}\label{fiblefcomp}
 If  $\con{h}{G}$ is
  set of elements of $G$ conjugate to $h$
  \[(\widehat{L_{X/G}})(f\cdot g)(\gamma)=\left\lbrace\begin{array}{ll}|\cen{h}{G}|&\text{if }\Theta(\gamma)=h\text{ and } g\in \con{h}{G}\\
  0&\text{otherwise}\end{array}\right.
  \]
  If $x\in X/G$ is a fixed point of $\bar{f}$ and $\tilde{x}$ is a lift of $x$ to $X$ so that $f(\tilde{x})=\tilde{x}h$ then 
  \[(\widehat{R_{X/G}})(f\cdot g)(c_x)=\left\lbrace\begin{array}{ll}\displaystyle\sum_{\{k\in G|k=hkg \}} i (\tilde{x}k)&\text{if } g\in \con{h}{G}\\
  0&\text{otherwise}\end{array}\right.
  \]
  where $i$ is the inclusion of the fixed points of $f\cdot g$ into  $\Lambda^{f\cdot g} X$ as constant paths.
\end{prop}

\begin{proof}
  For a $\bar{f}$ twisted loop $\gamma$ in $X/G$ and a lift $\tilde{\gamma}$ of $\gamma$ that satisfies $f(\tilde{\gamma}(0))=\tilde{\gamma}(1)h$ the image of $\tilde{\gamma}(1)$
  under the composite
  \begin{equation}
  F_{\bar{f}(x)}\to F_x\xto{ f\cdot g}F_{\bar{f}(x)}
  \end{equation}
  is $\tilde{\gamma}(1)hg$
  and the image of $\tilde{\gamma}(1)k$ is $\tilde{\gamma}(1)hkg$.  

Since $F_{\bar{f}(x)}$ is discrete the Lefschetz number is the number of points fixed 
by the endomorphism.    The group action of $G$ on $F_x$ is free, so we have a fixed point for each $k\in G$ where $k=hkg$.
  For any $k\in G$ of this form the map $l\mapsto lk^{-1}$ defines a  bijection from the centralizer of $h$ in $G$ 
  to 
  $\{l\in G|l=hlg\}$.

The Reidemeister trace is the sum of the constant paths associated to the fixed points.
\end{proof}

\section{Lefschetz numbers for spaces with free actions}\label{leffree}

As we see from the previous section, working with parametrized spaces is convenient and powerful.  Unfortunately, 
this approach does not immediately translate to invariants defined using more classical approaches such as \cite{FadellWong,Weber, Weber2, Wong, Wong2}.  
To compare the invariants defined here with these alternatives, we will follow the approach in \cite{KW} and 
replace parametrized spaces by spaces with a group action.   

To prove Theorem B for the Lefschetz number we start with a description of $L(f\cdot g)$ in terms of $R(\bar{f})$.  Let $\nu_g$ be the composite 
\[\Lambda^{\bar{f}}(X/G)_+\xto{\Theta} \sh{G}_+\xto{\mu_g} S^0.\]

\begin{lem}\label{infgeoorbits}
  For each $g\in G$, the stable map
  \[S^0\xto{R(\bar{f})} \sh{(X/G)_{\bar{f}}}\cong \Lambda^{\bar{f}}(X/G)_+\xto{\nu_g} S^0\]
  is $\frac{1}{|\cen{g}{G}|}(L(f\cdot g))$.
\end{lem}

\begin{proof}
  By \autoref{fiblefcomp} $\widehat{L_{X/G}}(f\cdot g)$ agrees with the composite 
  \[\sh{(X/G)_{\bar{f}}}\xto{\nu_g} S^0\xto{|\cen{g}{G}|} S^0.\]

   We have a commutative diagram \[\xymatrix{
  S^0\ar[r]^-{R(\bar{f})}\ar[dr]_{L(f\cdot g)}&\sh{(X/G)_{\bar{f}}}\ar[d]^{\widehat{L_{X/G}}(f\cdot g)}\\
  &S^0
  }\] as in \eqref{fib_reidemeister}
  and so $
  L(f\cdot g)=\widehat{L_{X/G}}(f\cdot g)\circ R(\bar{f})=|\cen{g}{G}|\nu_g\circ R(\bar{f}).$
\end{proof}

To identify $\mu_g\circ \tr_{/G}(f)$ and $\nu_g\circ R(\bar{f})$ we first provide another description of the $G$-space $X$.  Since $X\to X/G$ is a covering space, there is a classifying map $\phi\colon X/G\to BG$ 
and a (homotopy) pullback diagram
\[\label{pullbackdia}\xymatrix{
X\ar[r]^\psi\ar[d]^\pi&EG\ar[d]^\pi\\
X/G\ar[r]^\phi&BG
}\]
We let $(EG,\pi)$ be $EG\amalg BG$ regarded as a space over $BG\times \ast$ using $\pi$.  In 
the same way $(X,\pi)$ is $X\amalg X/G$ regarded as a space over $X/G\times \ast$.  
The diagram above defines  an equivariant fiberwise equivalence \cite[3.3]{PS_mult}
\[(X,\pi)\cong {}_\phi(BG)\odot (EG,\pi)\] 
and  
$X_+\cong S^0_{X/G}\odot {}_\phi(BG)\odot (EG,\pi)$.

\begin{lem}\label{factor}
  There is a fiberwise map $H\colon (X/G)_{\bar{f}}\odot {}_\phi(BG)\to {}_\phi (BG)$ and the composite 
  \begin{align*}S^0_{X/G}\odot {}_\phi (BG) &\odot (EG,\pi)\\
   &\xto{\bar{f}\odot \id\odot \id} 
   S^0_{X/G}\odot (X/G)_{\bar{f}}\odot{}_\phi (BG)\odot (EG,\pi)\\
   & \xto{\id\odot H\odot \id}S^0_{X/G}\odot {}_\phi (BG)\odot (EG,\pi)\\
   & \xto{\id\odot \id\odot (-)g}S^0_{X/G}\odot {}_\phi (BG)\odot (EG,\pi),
  \end{align*}
  is $f\cdot g\colon X_+\to X_+$.
\end{lem}

The fiberwise map $(-)g\colon (EG,\pi)\to (EG,\pi)$ is multiplication by $g$. 

\begin{proof}
  If we give $EG\times EG$ the diagonal $G$-action, the projection maps $EG\times EG\to EG$ are 
  $G$-equivariant and so there is a $G$-homotopy $K\colon EG\times EG\times I\to EG$ between the projections \cite[14.4.4]{tomDieck}.
  In the composite 
  \[\xymatrix@C=35pt{X\times I \ar[r]^-{(f \times \id)\times \id}&
  X\times X\times I\ar[r]^-{\psi\times \psi\times \id}&
  EG\times EG\times I\ar[r]^-K&EG}\]
  all maps are equivariant, so there is an induced map 
  $H\colon X/G\times I \to BG$ that is a homotopy from $\phi \circ \bar{f}$ to $\phi$.

  This homotopy defines a map 
  $H\colon S_{\bar{f}}\odot {}_\phi BG\to {}_\phi BG$ as in \cite[3.3]{PS_mult}.  The result then follows by explicit computation.
\end{proof}

\autoref{paradual} implies  $S^0_{X/G}$ is right dualizable and \autoref{basechangedual} implies  $_\phi BG$ is right dualizable.  
To use \autoref{factor} and \autoref{thm:compose-duality} to factor the trace 
 we need to describe the dual of $(EG,\pi)$.
If $G$ is a finite group $(EG,\pi)$ is right dualizable \cite[15.1.1]{MS} but this is not the dual we need here.

We let $\widecheck{(EG,\pi)}$  denote the fiberwise space $EG\amalg BG$ regarded as a space over $\ast\times BG$ via $\pi$. 
Note that $(EG,\pi)$ and $\widecheck{(EG,\pi)}$ both have $G$ actions.  We regard 
$(EG,\pi)$ as a space with a right $G$ action and $\widecheck{(EG,\pi)}$ as a space with a left $G$ action.

\begin{lem}\label{funnydual} There is a map   $\triangle_!S^0_{BG}\coloneqq(BG^I, ev_0\times ev_1)\to(EG,\pi)\wedge_G \widecheck{(EG,\pi)} $ over $BG\times BG$ and
  and a  $G\times G$-equivariant map
  $\widecheck{(EG,\pi)}\odot (EG,\pi) \to G_+$ so that the composites as in the definition of a dual pair are homotopic to identity maps
  through homotopies that respect both the group action and the parametrized structure.
\end{lem}

We use $\wedge_G$ to indicate quotienting by the diagonal group action after the smash product and 
$\odot $ to indicate the bicategory composition in $\Ex$.

While these do not define a dual pair in any of the bicategories in \autoref{bicateg}, the given structure  will allow us to use them in much the same way.

\begin{proof}
  Define a map $\triangle_!S^0_{BG}\to(EG,\pi)\wedge_G\widecheck{(EG,\pi)}$ by taking the quotient of the evaluation map 
  $ev_0\times ev_1\colon EG^I\to EG\times EG$ by the pointwise action of $G$ on the domain and the diagonal action on the codomain. 
  The map 
  $EG\times_{BG} EG\to G$ that takes a pair $(x,y)$ in $EG\times EG $ to the element $g\in G$ so that $xg=y$ defines a map $\widecheck{(EG,\pi)}\odot (EG,\pi) \to G_+$.

  The image of $(e,\gamma)$ under the composite 
  \[\xymatrix@C=15pt{\widecheck{(EG,\pi)}\odot \triangle_!S^0_{BG}\ar[r]& \widecheck{(EG,\pi)}\odot(EG,\pi)\odot\widecheck{(EG,\pi)}\ar[r]
  &G_+\wedge_G \widecheck{(EG,\pi)}
  }\]
  is $(g ,\tilde{\gamma}(1))$ where $\tilde{\gamma}$ is a lift of $\gamma$ to $EG$ and $g$ is the element of $G$ that takes $e$ to $\tilde{\gamma}(0)$.
  In particular, if we compose with the isomorphisms $\widecheck{(EG,\pi)}\cong \widecheck{(EG,\pi)}\odot \triangle_!S^0_{BG}$ and  $G_+\wedge_G \widecheck{(EG,\pi)}
  \cong \widecheck{(EG, \pi)}$
  the map is homotopic to the identity map.  The other composite is similar.
\end{proof}

\begin{proof}[Proof of Theorem B (Free Case, Lefschetz Number)]
  As observed above, the maps in \autoref{funnydual} are not a coevaluation and evaluation for a dual pair, but we can use \autoref{thm:compose-duality} and these maps to define 
  a coevaluation and evaluation
  \[S^n\to X_+\wedge_G DX\quad \text{and}\quad DX\wedge X_+\to G_+\wedge S^n\]
  as in \autoref{adamsdual}.  Since the trace is independent of the choice of dual we can equally well use this dual to compute $\tr_{/G}(f)$.

  Applying \autoref{thm:compose-duality} to the decomposition in \autoref{factor}, we see that  
   $\tr_{/G}(f)$  is the composite  
  \[\xymatrix{S^0\ar[r]^-{R(\bar{f})}&\sh{(X/G)_{\bar{f}}}\ar[r]^-{\tr(H)}&\sh{\triangle_! S^0_{BG}}\ar[r]^-{``\chi\textrm{''}}&\sh{G_+}}\]
  where the rightmost map is defined using the ``dual pair'' in \autoref{funnydual}.  
  By \cite[5.7]{PS_mult} the image of a twisted loop $\gamma\colon x\mapsto \bar{f}(x)$ in $X/G$ under $\tr(H)$ is $\phi(\gamma)\cdot \alpha(x)$ where
  $\alpha$ is the path from $\phi(\bar{f}(x))$ to $\phi(x)$ defined by the homotopy  in \autoref{factor}.
  The map $``\chi\textrm{''}$ lifts this path to $EG$ and assigns the value $g$ if multiplying the initial point of this lift 
  by $g$ is the terminal point of the lift.
  We obtain the same element by  first lifting $\gamma$ to $X$ and a assigning a group element  in the same way.
  Then $\mu_g\circ ``\chi\textrm{''} \circ \tr(H)=\nu_g$ and 
  \[\mu_g\circ \tr_{/G}(f)= \mu_g\circ ``\chi\textrm{''} \circ \tr(H)\circ R(\bar{f})=\nu_g\circ R(\bar{f}).\qedhere\]
\end{proof}

\section{Reidemeister traces for spaces with free actions}\label{reifree}
  
We can use very similar ideas to prove the corresponding result for the Reidemeister trace. 

\begin{lem}
  For each $g\in G$, the composite 
  \[S^0\xto{R(\bar{f})} \sh{(X/G)_{\bar{f}}}
  \xto{\xi_g}\Theta^{-1}(g)_+\]
  is $\frac{1}{\lvert \cen{g}{G}\rvert }\zeta_g(R(f\cdot g))$.
\end{lem}

\begin{proof}
  From  \autoref{fiblefcomp} if $\tilde{x}$ is a lift of $x$ so that  $f(\tilde{x})=\tilde{x}h$  
  \[(\widehat{R_{X/G}})(f\cdot g)(c_x)=\left\lbrace\begin{array}{ll}\displaystyle\sum_{k\in \cen{h}{G}}i_x(\tilde{x}k)&\text{if }g\in \con{h}{G}\\
  0&\text{otherwise}\end{array}\right. 
  \]
 
  When we compose with $\zeta_h$ we disregard paths whose value  under $\Theta$ is not $h$ 
and we do not distinguish between the paths where $\Theta(\gamma)=h$. In particular, 
  \begin{align}
  \zeta_g\circ (\widehat{R_{X/G}})(f\cdot g)(c_x)=\left\lbrace\begin{array}{ll}|\cen{h}{G}|c_{x}&\text{if }g\in \con{h}{G}\\
  0&\text{otherwise}\end{array}\right.
  \end{align}
  Similarly, $\xi_h$ is the identity on paths whose value  under $\Theta$ is  $h$ and we disregard the other paths.  Then  
	$|\cen{h}{G}|\xi_g(c_x)=\zeta_g\circ (\widehat{R_{X/G}})(f\cdot g) (c_x)$
  and 
  \[|\cen{h}{G}|\xi_g\circ R(\bar{f})=\zeta_g\circ (\widehat{R_{X/G}})(f\cdot g)\circ R(\bar{f})=\zeta_g\circ R(f\cdot g).\qedhere\] 
\end{proof}

To complete the proof of Theorem B it only remains to compare $\xi_g\circ \tr_{/G}(f)$ and $\xi_g \circ R(\bar{f})$.  This is 
essentially identical to the proof of the corresponding comparison in the previous section. 

\begin{proof}[Proof of Theorem B (Free Case, Reidemeister Trace)]
  The quotient  map $\tilde{X}\to X\to X/G$ is classified by a map 
  $\phi\colon X/G\to B\pi_1(X/G)$ and so we can write $\tilde{X}$ as the composite 
  \[S^0_{X/G}\odot {}_\phi (B\pi_1(X/G))\odot (E\pi_1(X/G),\pi)\] and $\tilde{f}\cdot g$
  as a composite
  \begin{align}
  S^0_{X/G}&\odot {}_\phi (B\pi_1(X/G))\odot (E\pi_1(X/G),\pi)\\
  &\xto{\bar{f}\odot \id\odot \id} S^0_{X/G}\odot (X/G)_{\bar{f}}\odot {}_\phi (B\pi_1(X/G))\odot (E\pi_1(X/G),\pi)\\
  &\xto{\id\odot H\odot \id} S^0_{X/G}\odot {}_\phi (B\pi_1(X/G))\odot (B\pi_1(X/G))_{B(f)}\odot (E\pi_1(X/G),\pi)\\
  &\xto{\id\odot \id\odot (-)g} S^0_{X/G}\odot {}_\phi (B\pi_1(X/G))\odot (E\pi_1(X/G),\pi)\wedge \pi_1(X/G)_{\bar{f}_*+}
  \end{align}
  The map $H\colon (X/G)_{\bar{f}}\odot {}_\phi (B\pi_1(X/G))\to {}_\phi (B\pi_1(X/G))\odot (B\pi_1(X/G))_{B(f)}$
  is induced by the square
  \[\xymatrix{
  X/G\ar[r]^-\phi\ar[d]^{\bar{f}}&B\pi_1(X/G)\ar[d]^{B(f)}\\
  X/G\ar[r]^-\phi&B\pi_1(X/G)
  }\] as in \cite[3.3]{PS_mult}.
  If $\gamma$ is a path in $B\pi_1(X/G)$ so that $\gamma(1)=\pi(e)$ then
  \[(-)g\colon  (B\pi_1(X/G))_{B(f)}\odot (E\pi_1(X/G),\pi)\to  (E\pi_1(X/G),\pi)\] is defined by taking 
  a pair $(\gamma,e)$ to $\gamma_*(f(e)g)$ where $\gamma_*$ is the map induced on fibers by $\gamma$.
  This map does not respect the $\pi_1(X/G)$ action unless we twist the action on the target by $f$. 

  The spaces $S_{X/G}^0$ and $_\phi (B\pi_1(X/G))$  are right dualizable, and $(E\pi_1(X/G),\pi)$ has a dual as in \autoref{funnydual}.
  This gives a factorization of $\tr_{/\pi_1(X/G)}(\tilde{f})$ as the composite
  \[S^0\xto{R(\bar{f})} \sh{X/G_{\bar{f}}}\xto{\tr(H)} \sh{(B\pi_1(X/G))_{B(f)}}\xto{``\chi\textrm{''}} \sh{\pi_1(X/G)_{f_*}}_+\]
  where the second and third maps take a twisted loop to its homotopy class.
  Then 
  \[ \sh{X/G_{\bar{f}}}\xto{\tr(H)} \sh{(B\pi_1(X/G))_{B(f)}}\xto{``\chi\textrm{''}} \sh{\pi_1(X/G)_{f_*}}_+\xto{\xi_g}\Theta^{-1}(g)_+\]
  is $\sh{X/G_{\bar{f}}}\xto{\xi_g} \Theta^{-1}(g)_+$.
\end{proof}

\part{General actions}

We now generalize to group actions that are not necessarily free.  
Using Theorems \ref{ulrichadd} and  \ref{radd1}, to understand $L_G(f)$ and $R_G(f)$ it is enough to understand the invariants $L_G(f_{(H)})$ and $R_G(f_{(H)})$ for each subgroup 
$H$ of $G$.   With this observation, the following statement is the relevant 
form  of Theorem A.
\begin{thmA}[General case]\label{facgenrei} 
    For each subgroup $H$ of $G$,
  $\chi(\WH)L_G(f_{(H)})=\chi_G(G/H)L(f_H)$ and 
  there integers $a_\gamma$ so that 
   \[R_G(f_{(H)})=\sum a_\gamma (i_\gamma \circ \tr^\triangle_G(F_\gamma)),\,\, R(f_{(H)})=\sum a_\gamma (i_\gamma \circ \tr^\triangle(F_\gamma)),\]
\[R_{\WH}(f_{H})=\sum a_\gamma (i'_\gamma \circ \tr_{\WH}^\triangle(F'_\gamma))\text{ and } R(f_{H})=\sum a_\gamma (i'_\gamma \circ \tr^\triangle(F'_\gamma))\]
  where all sums are taken over $\pi_0(\Lambda^{f^H/\WH}(X^H/\WH))$.  
\end{thmA}
The equivariant Lefschetz number of the endomorphism of $X^{(H)}/X^{>(H)}$ induced by $f$ is denoted $L_G(f_{(H)})$ and    
 $L(f_H)$ is the Lefschetz number of the endomorphism of $X^H/X^{>H}$ induced by  $f$. 
Here$F_\gamma$ is the fiber of $X^{(H)}\to X^H/\WH$ over $\gamma(1)$, $i_\gamma$ is the inclusion 
as constant paths, 
  $F_\gamma'$ is the fiber of $X^H\to X^H/\WH$ over $\gamma(1)$ and $i_\gamma'$ is the inclusion of 
the fiber over $\gamma(1)$ as constant paths.

For Theorem  B we need to give a more explicit description of the relevant traces.
In \autoref{isodualind} we will see that for each subgroup $H$ of $G$ the space $X^H\cup CX^{>H}$, the mapping 
cone of the inclusion $X^{>H}\to X^H$, is dualizable in the bicategory $\gTop$ as a $\WH$ space. 
The trace of the induced map \[f^H\colon X^H\cup C X^{>H}\to X^H\cup C X^{>H}\] is a map
$\tr_{/\WH}(f^H)\colon S^0\to \sh{\WH}_+$.  As before, 
this trace is carrying too much information.  For each conjugacy class $g$ in $\sh{\WH}_+$ there is a map $\mu_{\WH,g}\colon \sh{\WH}_+\to S^0$.

In \S\ref{rgen} we will take a very similar trace of an endomorphism of the universal cover of $X^H$.  The result is a map $\tr_{/\pi_1(X^H)\rtimes \WH}(f_{H})\colon S^0\to \sh{(\pi_1(X^H)\rtimes \WH)_{f^H}}_+$ and there are corresponding 
quotient maps $\xi_{\WH,g}$ and $\zeta_{\WH,g}$.

\begin{thmB}[General case] \label{isotraceind}
  For each $g\in \WH$, \[L(f_H\cdot g) = \lvert \cen{g}{\WH}\rvert \mu_{\WH,g} (\tr_{/\WH}(f_H)) \]
  \[\zeta_{\WH,g}R(f_{H}\cdot g)=\lvert \cen{g}{\WH}\rvert \xi_{\WH,g} \tr_{/\pi_1(X^H)\rtimes \WH} (\tilde{f}_{H}).\] 
\end{thmB}

One significant difference between the free case and the general case  is that from this point on we will restrict to closed smooth manifolds or finite CW complexes.  This reflects
  a single result, \autoref{CWdualisotropy}, which may have a generalization for ENRs.

\section{Homotopical Invariants}\label{hoinvgen}

The proof of Theorem A follows the proof of  the corresponding free result but requires a refinement of the relative Reidemeister trace \cite{relative}.
For a map $Y\to B$ and $A\subset Y$, $C_B(Y,A)$ is the homotopy pushout of the maps 
\[Y\amalg B\leftarrow A\amalg A\to A.\]  This is regarded as a parametrized space over $\ast\times Y$.  
The {\bf relative Reidemeister trace} of a map $f\colon Y\to Y$ so that $f(A)\subset A$ is the bicategorical trace of the 
induced map \[C_Y(Y,A)\to C_Y(Y,A)\odot Y_f.\]

\begin{lem}\label{CWdualisotropy}
  Let $Y$ be a simplicial complex and $A\subset Y$ be a subcomplex.  There is a neighborhood $U$ of $A$ in $Y$ so that 
  $C_{Y\setminus A}(Y\setminus A, U\setminus A)$ is right dualizable.
\end{lem}

It is possible that a similar result would hold for ENRs, but the only proof I know of at this time makes significant use of the simplicial structure.

\begin{proof}
  Choose a mapping cylinder neighborhood $U$ of $A$ in $Y$ \cite{daverman, whitehead}.  
  There is a retraction of $ U\setminus A$ to $\partial (U)$ and a corresponding retraction 
  $Y\setminus A\to Y^0\coloneqq Y\setminus \mathrm{Int}(U)$.  This gives an equivalence \[(Y\setminus A, U\setminus A)\cong (Y^0, \partial U).\] 
  Since $(Y^0, \partial U)$ is a compact CW pair it is a compact ENR pair and \cite[18.5.2]{MS} implies 
  that $C_{Y^0}(Y^0, \partial U)$ is right dualizable.  

  If $i\colon Y^0\to Y\setminus A$ is the inclusion $C_{Y^0}(Y^0, \partial U)\odot {}_i(Y\setminus A)\simeq C_{Y\setminus A}(Y^0, \partial U)$ \cite[p. 1288]{PS_mult}
  and \autoref{basechangedual} implies $C_{Y\setminus A}(Y^0, \partial U)$ 
    is right dualizable.  
By excision \cite[18.4.5]{MS}
  \[C_{Y\setminus A}(Y^0, \partial U)\simeq C_{Y\setminus A}(Y\setminus A, U\setminus A).\qedhere\]
\end{proof}

\begin{lem}\label{relfactor}
  If $A\subset Y$ are closed smooth manifolds and $f\colon Y\to Y$ is a continuous map so that $f(A)\subset A$ 
  there is a class in $\pi_0^s(\Lambda^f (Y\setminus A))$ so that the image in $\pi_0^s(\Lambda ^fY)$ is the relative Reidemeister trace of $f$.
\end{lem}

Following the notation above $\Lambda^f (Y\setminus A)\coloneqq\{(\gamma,y)\in (Y\setminus A)^I\times Y|\gamma(1)=f(y)\}$. From a classical view of the Reidemeister trace in 
terms of fixed point indices and fixed point classes this is a very intuitive statement since 
the constant paths at the fixed points in $Y\setminus A$
are elements of $\pi_0(\Lambda^f (Y\setminus A))$ as well as elements of $\pi_0(\Lambda ^fY)$.  

\begin{proof}
  If $i\colon Y\setminus A \to Y$ is the inclusion, excision and \cite[18.4.4]{MS} imply \[C_{Y}(Y,A)\simeq C_{Y\setminus A}(Y\setminus A, U\setminus A)\odot {}_i Y.\]
  The induced map 
  $f\colon C_{Y}(Y,A)\to C_{Y}(Y,A)\odot Y_{f}$ defines a map 
  \[\xymatrix@R=1pt{C_{Y\setminus A}(Y\setminus A, U\setminus A)\odot {}_i Y\ar[r]& C_{Y\setminus A}(Y\setminus A, U\setminus A)\odot {}_i Y\odot  Y_{f}.}\] 
  Via adjunction and the dual pair in \autoref{basechangedual}, we have a map 
  \begin{equation}
  \xymatrix@R=1pt{C_{Y\setminus A}(Y\setminus A, U\setminus A)\ar[r]& C_{Y\setminus A}(Y\setminus A, U\setminus A)\odot {}_i Y\odot Y_{f\circ i}.
  }\label{genreldecom}\end{equation}
  \autoref{thm:compose-duality}  implies
  the diagram below where the vertical map is the inclusion of paths and $R_A(f)$ is the relative Reidemeister trace of $f$ with respect to $A$ commutes.
  \[\xymatrix{
  S^0\ar[r]^-{\tr\eqref{genreldecom}}\ar[dr]_{R_A(f)}&\sh{ {}_i Y_{f\circ i}}\ar[d]\\
  &\sh{ Y_{f}}
  }\] 

  We can lift $R_A(f)$ further.  Since $A$ is an NDR in $Y$ there is a  map $u\colon Y\to I$ so that $u^{-1}(0)=A$ and $u^{-1}([0,1))=U$ for some 
neighborhood $U$ of $A$ in $Y$.  Then  one choice for  \eqref{genreldecom} is
  \[(y,t)\mapsto ((f(x),tu(f(y)) ),(f(x),c_{f(x)},x)). \]
  Note that if $y\in Y$ satisfies  $f(y)\in A$ then the image of $(y,t)$ is in the section.  In fact, the only points whose images are not in the section are those 
    for fixed points in $Y\setminus A$.  If $P$ is the subspace
  $(Y\setminus A)\times_{i}(Y^I)\times_f (Y\setminus A)$ consisting of triples where the path lies entirely in $Y\setminus A$ and $P$ is regarded as a space over $(Y\setminus A)\times (Y\setminus A)$
  \eqref{genreldecom}  lifts to define a map 
  \[C_{Y\setminus A}(Y\setminus A, U\setminus A)\to C_{Y\setminus A}(Y\setminus A, U\setminus A)\odot P_+ \] where $P_+$ is $P$ with a disjoint section added.  
  This allows us to factor  the Reidemeister trace of $f$ through the twisted loops in $Y\setminus A$.
\end{proof}

Let   $i^H\colon X_H\coloneqq\{x\in X\mid G_x=H\}\to X^H$ be the inclusion and $\pi_H\colon X_H\to X_H/\WH$ be the quotient map.
If $U^H$ is an open neighborhood of $X^{>H}$ in $X^H$ that 
retracts to $X^{>H}$   and \[C_HX\coloneqq C_{X_H/\WH}(X_H/\WH, (U^H-X^{>H})/\WH)\]
 excision implies  \cite[18.4.5]{MS}    $C_{X^H}(X^H, X^{>H})$ can be written as the composite
\begin{equation}\label{genreldecom2}C_HX\odot (X_H/\WH)_{\pi_H}\odot {}_{i^H}(X^H).\end{equation}

\begin{lem}\label{facgenrei2}
  The map $C_{X^H}(X^H, X^{>H})\to C_{X^H}(X^H, X^{>H})\odot (X^H)_{f^H}$ induced by $f^H$ factors as a composite 
  \begin{align}C_HX\odot (X_H/&\WH)_{\pi_H}\odot {}_{i^H}(X^H)\to C_HX\odot P_+\odot (X_H/\WH)_{\pi_H}\odot {}_{i^H}(X^H)\\
  &\to C_HX\odot (X_H/\WH)_{\pi_H}\odot {}_{i^H}(X^H)\odot (X^H)_{f^H}\odot (X^H)_{i^H}\odot {}_{i^H}(X^H)\\
  &\to C_HX\odot (X_H/\WH)_{\pi_H}\odot {}_{i^H}(X^H)\odot (X^H)_{f^H}
  \end{align}
\end{lem}
Following the notation above $P\coloneqq X_H\times_{i^H} (X_H)^I\times_{f^H} X_H$.  

\begin{proof}  The first map is as in the proof of \autoref{relfactor} and 
  the last map is the evaluation for the dual pair in \autoref{basechangedual}. 

  If $P(f,g)$ is the homotopy pullback of maps $f$ and $g$ we can define a map 
  \[\xymatrix@R=1pt{P\times_{X_H/\WH}P(\id,\pi_H)\times_{X_H}P({i^H},\id)\ar[r]& P(\id,\pi_H)\times_{X_H}P({i^H},\id)\times_{X^H} P(\id,f^H)\\
  (\gamma,x),(\beta,y),(y, \alpha)\ar@{|->}[r]& (\gamma, f(\tilde{\beta}(0))),(f(\tilde{\beta}(0)), f(\tilde{\beta})),(f(\alpha),\alpha(1))}
  \] where
  $\tilde{\beta}$ is the lift of $\beta$ to a path in $X_H$ so that $\tilde{\beta}(1)=y$. 
  This induces a map 
  \[P_+\odot (X_H/\WH)_{\pi_H}\odot{}_{i^H}(X^H)\to (X_H/\WH)_{\pi_H}\odot {}_{i^H}(X^H)\odot (X^H)_{f^H}.\]

  Then we have a map 
  \begin{align}
  P_+\odot (X_H/\WH)_{\pi_H}&\xto{\sim} P_+\odot (X_H/\WH)_{\pi_H}\odot \triangle_!S^0_{X_H/\WH}\\
  &\to P_+\odot (X_H/\WH)_{\pi_H}\odot{}_{i^H}(X^H)\odot (X^H)_{i^H}\\
  &\to  (X_H/\WH)_{\pi_H}\odot {}_{i^H}(X^H)\odot (X^H)_{f^H}\odot (X^H)_{i^H}
  \end{align}
  where the first map is an isomorphism, the second is the coevaluation for the dual pair in \autoref{basechangedual}, and the third is defined above.
\end{proof}

Note that the same argument applies to the map \[C_{X^{(H)}}(X^{(H)}, X^{>(H)})\to C_{X^{(H)}}(X^{(H)}, X^{>(H)})\odot (X^{(H)})_{f^{(H)}}.\]
We will also be able to use identical arguments to prove both version of the Reidemeister trace statements in Theorem A.  Because of the 
remarkable similarity we only describe the first case.

\begin{proof}[Proof of Theorem A (General Case)]
  Using \autoref{CWdualisotropy}
   we can choose a neighborhood $U^{H}$ of $X^{>H}$ in $X^{H}$ so that $C_HX$ is dualizable. 
  \autoref{basechangedual} implies  ${}_{i^{(H)}}X^{(H)}$ is right dualizable.
  Since $\pi_{(H)}\colon X_{(H)}\to X_H/\WH$ is a fibration with finite fiber $(X_H/\WH)_{\pi_{(H)}}$ is right dualizable by   \cite[4.7]{PS_mult}. 
  Then $C_{X^{(H)}}(X^{(H)}, X^{>(H)})$ is right dualizable and \autoref{thm:compose-duality} and \autoref{facgenrei2} give us a decomposition of the relative Reidemeister trace as the composite 
 \[S^0\to \sh{P}\to \sh{{}_{i^{(H)}}X^{(H)}\odot X^{(H)}_{f^{(H)}}\odot X^{(H)}_{i^{(H)}}}\to \sh{X^{(H)}_{f^{(H)}}}\]
  where the first map is the lift of the relative Reidemeister trace of $\bar{f}^{(H)}$ to the twisted loops in $X_H/\WH$. 
  The second associates to a twisted loop in $X_H/\WH$ the Reidemeister trace of the induced endomorphism of the fiber.  (There is an induced endomorphism
  since the path is entirely contained in $X_H/\WH$ and $X_{(H)}\to X_H/\WH$ is a fibration.)  The third map is the inclusion of paths.

  The computation of the Reidemeister trace of the fiber is as in the proof of \autoref{fibcompute}.  Since the factoring applies in both the
  equivariant and nonequivariant cases we have the desired decompositions of the Reidemeister trace.

  If we choose the integers $a_x$ so that $R(f_H/\WH)=\displaystyle\sum_{x\in \Fix(\overline{f}_H)} a_x c_x$ is the relative Reidemeister trace of $f_H/\WH$ relative to the subspace $X^{>H}/\WH$
  we have the decomposition in the statement of the theorem.

  If we compose with the collapse maps $\Lambda^{f^H}X^H\to \ast$, identifications of the Reidemeister trace gives corresponding identifications of the Lefschetz number:
  \begin{align}L_G(f_{(H)})&=\sum_{x\in \Fix(\overline{f_H})(e)} a_x \chi_G(F_x)=\chi_G(G/H) \sum_{x\in \Fix(\overline{f_H})(e)} a_x  \\
  L(f_{(H)})&=\sum_{x\in \Fix(\overline{f_H})(e)} a_x  \chi(F_x)=\chi(G/H)\sum_{x\in \Fix(\overline{f_H})(e)} a_x\end{align}
  and $\chi(G/H)L_G(f_{(H)})=\chi_G(G/H)L(f_{(H)})$.

  To replace $L(f_{(H)})$ by a multiple of $
  L(f_H)$ note that \[X^{(H)}/X^{>(H)}\cong \displaystyle\bigvee_{K\in \con{H}{G}} X^K/X^{>K}\text{ and so }
  L(f_{(H)})=\displaystyle\sum_{K\in \con{H}{G}} L(f_K).\] Traces are invariant under cyclic permutation, so if $K$ and $H$ are conjugate $L(f_K)=L(f_H)$.  Then 
  $\chi(\WH) L(f_{(H)})= \chi(G/H)L(f_H)$
  and so 
  \[\chi(\WH)\chi(G/H)L_G(f_{(H)})=\chi(\WH)\chi_G(G/H)L(f_{(H)})=\chi(G/H)\chi_G(G/H)L(f_H).\]
  \end{proof}

\section{Lefschetz numbers for $G$-spaces}\label{lefgen}
Much of the work required to prove the general case of Theorem B has been done in the previous section 
and in the corresponding results in Part 1. 

\begin{lem}\label{isodualind}
  If $X$ is a compact $G$-ENR or closed smooth $G$-manifold then $X^H\cup CX^{>H}$ is dualizable as a $\WH$-space in $\gTop$.
\end{lem}

It is important to note that we cannot expect dualizability for $X^H$ in $\gTop$ if the action of $H$ is not free.   This is analogous to 
the algebraic requirement that dualizable modules must be projective.  

\begin{proof}
  If $X$ is a compact $G$-ENR  then $X^H$ is a closed $\WH$-ENR.   As in the proof of  \autoref{ulrichadd},
  $X^{>H}$ is a $\WH$-ENR and the inclusion $X^{>H}\to X^H$ is a $\WH$-cofibration.    
  Using \autoref{smcadd} $X^H/X^{>H}$ is $\WH$-dualizable and since the action of $\WH$ is free away from the basepoint 
  \cite[8.6]{Adams} implies that $X^H/X^{>H}$ is dualizable in $\gTop$ and the dual agrees with the $\WH$-equivariant dual.
\end{proof}

This statement is essentially the only change we need to make in the proof of the free case of Theorem B for the Lefschetz number.
 
\begin{proof}[Proof of Theorem B (General case, Lefschetz number)]

  If $U^H$ is a neighborhood of $X^{>H}$ in $X^H$ that retracts onto $X^{>H}$
  and $U_H\coloneqq U^H\setminus X^{>H}$  then \[X^H\cup C(X^{>H})\simeq X_H\cup C(U_H).\] 
  Since the action of $\WH$ on $X_H$ is free, there is a map $\phi\colon X_H/\WH\to B\WH$ that 
  classifies $X_H\to X_H/\WH$.  We can choose $U_H$ so that it is $\WH$ equivariant and it is classified by the restriction of 
  $\phi$ to $U_H/\WH$.  The result then follows from the proof of the free case if we replace $G$ by $\WH$ and $S^0_{X}$ by $C_{X_H/\WH}(X_H/\WH, U_H/\WH)$. 
\end{proof}  
  
We can collect these traces $\tr_{/\WH}(f_H)$ into a single trace by generalizing from the category $\gTop$ to the  category of \emph{profunctors}.  
Associated to a symmetric monoidal category $\sV$ with unit $S$ and monoidal product $\otimes$
there is a bicategory  $Pro(\sV)$  where
\begin{itemize}
  \item The objects are small categories,
  \item Between two small categories $A$ and $B$ we have the category of functors \[A\times B\op\to \sV\]
  and their natural transformations,
  \item For any small category $A$, there is a functor $U_A\colon A\times A\op\to \sV$ defined by 
  $U_A(a,a')=\amalg_{A(a,a')}S$.
  \item For functors $X\colon A\times B\op\to \sV$ and $Y\colon B\times C\op\to \sV$, we define $X\odot Y\colon A\times C\op\to\sV$ by taking 
   $(X\odot Y)(a,c)$  to be the coequalizer of the diagram
  \[\xymatrix{\displaystyle{\coprod_{b\to b'\in B}} X(a,b)\otimes Y(b',c)\ar@<.5ex>[r]\ar@<-.5ex>[r]&
  \displaystyle{\coprod_{b\in\ob \sB}}X(a,b)\otimes Y(b,c)}\] where the maps are induced by the action of the morphisms of $B$ on $X$ and $Y$.
  \item The shadow of $X\colon A\times A\op\to \sV$ is the coequalizer of the diagram \[\xymatrix{\displaystyle{\coprod_{a\to a'\in \mathrm{ob}(\sA)}}Z(a,a')
  \ar@<.5ex>[r]\ar@<1.5ex>[r]&\displaystyle{\coprod_{a\in \mathrm{ob}(A)}} Z(a,a)}.\]
  The symmetry isomorphism in $\sV$ defines the map 
  \[\sh{X\odot Y}\rightarrow \sh{Y\odot X}.\]
\end{itemize}

For topological examples we replace coequalizers by homotopy coequalizers using  the  bar resolution.  As in $\Top$, $\gTop$, and $\Ex$, 
we say a functor $A\times B\op\to \Top$ is dualizable if the composite with the suspension spectrum functor is dualizable.  This can 
be described using natural transformations 
$\eta\colon S^n\wedge A_+\rightarrow X\odot Y$ and 
$\epsilon\colon Y\odot  X\rightarrow S^n\wedge B_+$ in $\Top$ so that the usual duality composites are the identity after suspension by 
a sufficiently large sphere.

The generalization of the group $G$ is the component category. 
\begin{defn}
  The {\bf equivariant component category} $\ecomcat{G}{X}$ for a  $G$-space $X$ has objects  $G$-maps $x(H)\colon G/H\rightarrow X$.  
  The morphisms from $x(H)$ to $y(K)$ are the  $G$-maps \[\alpha\colon
  G/H\rightarrow G/K\] such that $y(K)\circ \alpha$ and $x(H)$ are $G$-homotopic.
\end{defn} 
For any $g\in G$ so that $g^{-1}Hg\subset K$ there is a $G$-map $R_g\colon G/H\to G/K$
defined by $R_g(lH)=lgK$.   All $G$-maps $G/H\to G/K$ are of this form and two such maps $R_g$ and $R_h$ are the same only if 
$gh^{-1}\in K$.

  If $x(H)\colon G/H\rightarrow X$ is a $G$-map and $C_x$ is the component of $X^H$ that contains $x(eH)$ let
  $X^H(x)$ be the pullback of the quotient map $X^H\to X^H/\WH$ and the inclusion $(\WH C_x)/\WH\to X^H/\WH$. 
  If $X^{>H}(x)\coloneqq \{y\in X^H(x)|H \subsetneq G_y\}$ 
define a functor \[\ecomsp{G}{X}\colon \ecomcat{G}{X}\op\rightarrow \Top\] 
by $\ecomsp{G}{X}(x(H))=X^H(x)/X^{>H}(x)$.  On morphisms we use  the induced group action.  

\begin{prop}\label{comdual} If $X$ is a compact $G$-ENR or closed smooth $G$-manifold then $\ecomsp{G}{X}$ is dualizable.
\end{prop} 
\begin{proof}
Using  \cite[3.7]{relative} it is enough to show that $\ecomsp{G}{X}(x(H))$ is dualizable for each object $x(H)$ relative to 
the action by $\ecomcat{G}{X}(x(H),x(H))$.  This is \autoref{isodualind}.
\end{proof}

Given an equivariant map $f\colon X\rightarrow X$
let $\ecomcatf{G}{X}$ be the functor $\ecomcat{G}{X}\times \ecomcat{G}{X}\op\to {\Top}$ defined by 
\[\ecomcatf{G}{X}(x(H),y(K))\coloneqq \ecomcat{G}{X}(f(y(K)),x(H)).\]
Then an endomorphism $f\colon X\rightarrow X$ induces a natural transformation 
\begin{equation}\label{lgenmap2}\ecommap{f}\colon \ecomsp{G}{X}\rightarrow \ecomsp{G}{X}\odot \ecomcatf{G}{X}.\end{equation}
The trace of $\ecommap{f}$ is a map $S^0\to \sh{\ecomcatf{G}{X}}$.

\begin{prop}  If  $\WH_{x,f}\coloneqq\{g\in \WH| [f(x)g]=[x]\in \pi_0(X^H)\}$ and $B(X)$ 
is the isomorphism classes of objects of $\ecomcat{G}{X}$ there is an isomorphism \begin{equation}\delta\colon 
\sh{\ecomcatf{G}{X}}_+\rightarrow \displaystyle\coprod_{x(H)\in B(X)} \sh{\WH_{x,f}}_+.
\end{equation}
  If $\delta_{x(H),f}\colon \sh{\ecomcatf{G}{X}}_+\rightarrow \sh{\WH_{x,f}}_+$ is one of the projections,  the composite 
  \[S^0\xto{\tr(\bar{f})} \sh{\ecomcatf{G}{X}}\xto{\delta_{x(H),f}} \sh{\WH_{x,f}}\] 
  is $\tr_{/\WH}(f)$. 
\end{prop}

\begin{proof}
Using the identifications 
\begin{align}
\ecomcat{G}{X}(f(x(H)),x(H))&=\{R_g\colon G/H\to G/H\mid f(x(H))\circ R_g\sim_G x(H) \}\\
&=\{g\in \WH\mid [f(x(H))g]=[x(H)]\in \pi_0(X^H)\}
\end{align} we have the isomorphism $\delta$ above.  Then 
\cite[3.6, 3.7]{relative} completes the proof.
\end{proof}

We can also compose with the cellular chain complex functor to define a functor
\[C_*(\ecomsp{G}{X})\colon \ecomcat{G}{X}\op\rightarrow \Ch_\mathbb{Q}.\]  If $X$ is a compact $G$-ENR or 
closed smooth $G$ manifold this functor is dualizable by functoriality.
The natural transformation $\ecommap{f}\colon \ecomsp{G}{X}\rightarrow \ecomsp{G}{X}\odot \ecomcatf{G}{X}$
induces a natural transformation 
\begin{equation}\label{lgenmap}
{\ecommap{f}}_*\colon C_*\left(\ecomsp{G}{X}\right)
\rightarrow C_*\left(\ecomsp{G}{X}\right)\odot \mathbb{Q}\ecomcatf{G}{X}\end{equation}
where $\mathbb{Q}\ecomcatf{G}{X}$ is  is the functor $\ecomcat{G}{X}\times\ecomcat{G}{X}\to \Ch_\mathbb{Q}$ 
defined by \[\mathbb{Q}\ecomcatf{G}{X}(x(H),y(K))\coloneqq \mathbb{Q}\ecomcat{G}{X}(f(y(K)),x(H)).\]
Functoriality of the trace implies the trace of \eqref{lgenmap} agrees with the trace of \eqref{lgenmap2}.  The
following result is an immediate consequence.

\begin{thm} If $f\colon X\rightarrow X$ is an equivariant map and the set
  \[\{x\in X|\text{there is }g\in G \text{ such that }f(x)=xg\}\]   is empty then the trace of  \eqref{lgenmap} is trivial.
\end{thm}

\section{Reidemeister traces for $G$-spaces}\label{rgen}

We can now combine \S\ref{reifree} and \S\ref{lefgen} and give another description of the Reidemeister trace for spaces with a group action that is not necessarily free.
We start with the common generalization of the fundamental group and the component category.
\begin{defn}
   The objects of the  {\bf equivariant fundamental category} $\efuncat{G}{X}$ of a $G$-space $X$ are the $G$-maps
  $x(H)\colon G/H\rightarrow X$.  A morphism from $x(H)$ to $y(K)$ is a $G$-map 
  \[R_g\colon G/H\rightarrow G/K\]
   and a homotopy class of $G$-maps 
  \[w(H)\colon G/H\times I \rightarrow X\] relative to 
  $G/H\times\partial I$ such that $w(H)(-,0)=x(H)$ and $w(H)(-,1)=y(K)\circ R_g$.

  The composite of $(R_g, w(H))$ and $(R_h, v(K))$ is $(R_h\circ R_g, (v(K)\circ R_g)w(H))$.
\end{defn}

Let $\widetilde{X^H}(x)$ be the universal cover of $X^H(x)$.   The usual action of paths on the cover defines an action of an endomorphism $(R_g, w(H))$  of $x(H)$ in $\efuncat{G}{X}$ on 
a point $\tilde{x}$ in 
$\widetilde{X^H}(x)$ by $\tilde{x}\mapsto (\tilde{x}g)\cdot w(H)(e)$.
Let $\rcov{X^{>H}}(x)$ be the pullback of $\widetilde{X^H}(x)$ along the 
inclusion $X^{>H}\to  X^H$, and  $\rcov{X_{H}}(x)$ be the pullback of $\widetilde{X^H}(x)$ along the 
inclusion $X_H\to  X^H$.  
The group action of 
$\efuncat{G}{X}(x(H),x(H))$ defines a group action on $\rcov{X_{H}}(x)$.
Define a functor \[\hat{X}\colon \efuncat{G}{X}\to \Top\] by $\hat{X}(x(H))\coloneqq\widetilde{X^H}(x)\cup C (\overline{X^{>H}}(x))$.

\begin{lem}\label{isodualr} If $X$ is a compact $G$-ENR or a closed smooth $G$-manifold then $\hat{X}(x(H))$
  is dualizable in $\gTop$ as a $\efuncat{G}{X}(x(H),x(H))$ space.
\end{lem}

\begin{proof}
  There is a diagram
  \[\xymatrix{{\rcov{X_H}(x)}\ar[r]^-{\tilde{\pi}_H}\ar[d]&X_H(x)\ar[r]^-{\pi_H}\ar[d]^{i^H}&X_H(x)/\WH\ar@{.>}[d]\\
   {\widetilde{X^H}(x)}\ar[r]^-{\tilde{\pi}^H}&X^H(x)\ar@{.>}[r]&X^H(x)/\WH
  }\]
  where the vertical maps are inclusions and the horizontal maps are quotients.  The top corner is a pullback square and hence a homotopy pullback since $\widetilde{X^H}\to X^H$ 
  is a covering map.  If $\pi\colon X_H(x)\to \ast$ then 
  \begin{align}(X_H(x)/\WH)_{\pi_H}\odot {}_{i^H}X^H(x) \odot (X^H(x))_{\tilde{\pi}^H}\odot {}_\pi\ast&\simeq (X_H(x)/\WH)_{\pi_H\tilde{\pi}_H}) \odot {}_\pi\ast\\&\simeq \widecheck{(\overline{X_H}(x),\pi_H\tilde{\pi}_H}_+. 
  \end{align}

  There is an isomorphism 
  \[\xymatrix@R=3pt{\efuncat{G}{X}(x(H),x(H))\ar[r]&\pi_1(X^H)\rtimes \WH\\
  (R_g,w(H))\ar[r]&(w(eH), g)
  }\] and this is compatible with the actions of each group on $\overline{X_{H}}(x)$. 
  The composite  of quotient maps
  $\overline{X_{H}}(x)\to X_{H}(x)\to X_{H}(x)/\WH$ is the quotient by the action of $\pi_1(X^H)\rtimes \WH$.
  Then we can use the approach of \autoref{funnydual} to define maps 
  \[\triangle_!S^0_{X_H/\WH}\to   \widecheck{(\overline{X_H}(x),\pi_H\tilde{\pi}_H)} \wedge_{\pi_1(X^H)\rtimes \WH}  (\overline{X_H}(x),\pi_H\tilde{\pi}_H)\]\[
  (\overline{X_H}(x),\pi_H\tilde{\pi}_H)  \odot\widecheck{(\overline{X_H}(x),\pi_H\tilde{\pi}_H)}\to (\pi_1(X^H)\rtimes \WH)_+\] 
  so that the required triangle diagrams commute.

  Using the decomposition before \autoref{facgenrei2} and {\cite[18.4.4]{MS}}
  \begin{align}
  \widetilde{X^{H}}(x)\cup &C(\rcov{X^{>H}}(x))\simeq C_{\widetilde{X^{H}}(x)}(\widetilde{X^{H}}(x),\rcov{X^{>H}}(x))\odot {}_\pi\ast\\
  &\simeq C_{X^H(x)}(X^H(x), X^{>H}(x)) \odot (X^H(x))_{\tilde{\pi}^H}\odot {}_\pi\ast\\
 &\simeq C_{H}X(x)\odot (X_H(x)/\WH)_{\pi_H}\odot {}_{i^H}X^H(x) \odot (X^H(x))_{\tilde{\pi}^H}\odot {}_\pi\ast
  \end{align} 
  From \autoref{CWdualisotropy} $C_{H}X(x)\coloneqq  C_{X_H(x)/\WH}(X_H(x)/\WH, (U^H-X^{>H})/\WH)$ is right dualizable.   The remaining parts of the 
  decomposition are dualizable by the discussion above.
\end{proof}

\autoref{isodualr} extends in the same way that \autoref{comdual} follows from \autoref{isodualind}.

\begin{prop}\label{geodual} 
  If $X$ is a compact $G$-ENR then $\hat{X}$ is dualizable as a right $\efuncat{G}{X}$-module.
\end{prop}

If $\efuncatf{G}{X}$ is the functor $\efuncat{G}{X}$-$\efuncat{G}{X}\to \Top$ defined by 
\[\efuncatf{G}{X}(x(H),y(K))\coloneqq \efuncat{G}{X}(f(y(K)),x(H)).\]
an equivariant map $f\colon X\rightarrow X$ defines a  natural transformation 
$\tilde{f}\colon \hat{X}\rightarrow \hat{X}\odot \efuncatf{G}{X}$.
The trace of $\tilde{f}$ is a map $\sh{S^0}\to \sh{\efuncatf{G}{X}}$.

The following result is a consequence  of \cite[3.6, 3.7]{relative}.

\begin{prop}  
There is an isomorphism $\sh{\efuncatf{G}{X}}\cong \coprod \sh{\efuncatf{G}{X}}(x(H),x(H))$ 
where the coproduct is taken over a choice of representatives of the isomorphism classes of objects of $\efuncat{G}{X}$.
The image of $\tr(\tilde{f})$ under the projection \[\coprod \sh{\efuncatf{G}{X}}(x(H),x(H)) \to \sh{\efuncatf{G}{X}}(x(H),x(H))\] is 
the trace in $\gTop$ of the induced map $\hat{X}(x(H))\to \hat{X}(x(H))\odot \efuncatf{G}{X}(x(H),x(H))$ with respect to the group action by $\efuncat{G}{X}(x(H),x(H))$.
\end{prop}

Exactly as before we can compose the dual pair for $\hat{X}$ with the cellular chain complex functor and define algebraic invariants.  This is
the  \emph{refined equivariant Lefschetz number} from \cite[5.7]{Weber}.
Using functoriality of the trace these agree with the topologically defined trace of $\tilde{f}\colon \hat{X}\rightarrow \hat{X}\odot \efuncatf{G}{X}$.

We can now finish the proof of Theorem B.  The map $\Theta$ defined in the Part 1 extends to a map \[\Theta_{\WH}\colon \sh{\efuncatf{G}{X}(x(H),x(H))}\to \sh{\WH}\]  by  $\Theta_{\WH}(R_g,w)= g$.   
Let \[\sh{\efuncatf{G}{X}(x(H),x(H))}_+\xto{\xi_{\WH,g}} \Theta_{\WH}^{-1}(g)_+\text{ and } \sh{\pi_1(X^H)_{gf^H_*}}_+\xto{\zeta_{\WH,g}} \Theta_{\WH}^{-1}(g)_+\] be collapse maps
generalizing the maps $\xi_g$ and $\zeta_g$.

\begin{proof}[Proof of Theorem B (General case, Reidemeister trace)]
  We can use essentially the same proof as in the case of a free action.

  The trace $\tr_{/\pi_1(X^H)\rtimes \WH}(\tilde{f}_{(H)})$ is the trace of the map 
  \[\widetilde{X^{H}}(x)\cup C(\rcov{X^{>H}}(x))\to \left(\widetilde{X^{H}}(x)\cup C(\rcov{X^{>H}}(x))\right)\wedge_{\pi_1(X^H)\rtimes \WH} (\pi_1(X^H)\rtimes \WH)_{f^H_*}
 \] induced by $f^H$.
Using the discussion above, $\widetilde{X^{H}}(x)\cup C(\rcov{X^{>H}}(x))$ is the composite
$C_{H}X(x)\odot (X_H(x)/\WH)_{\pi_H}\odot {}_{i^H}X^H(x) \odot (X^H(x))_{\tilde{\pi}^H}\odot {}_\pi\ast$.  We can replace this by $C_{H}X(x)\odot (X_H(x)/\WH)_{\pi_H}\odot  (X_H(x))_{\tilde{\pi}_H}\odot{}_{\tilde{i}^H}\overline{X^H}(x) \odot {}_\pi\ast$ where $\tilde{i}^H$
is the induced map $\overline{X^H}(x)\to\widetilde{X^H}(x)$.  We can further simplify 
to  $C_{H}X(x)\odot (X_H(x)/\WH)_{\pi_H\circ{\tilde{\pi}_H}}\odot {}_{\pi\circ\tilde{i}^H}\ast$.

There is a map $\phi\colon X_H(x)/\WH\to B(\pi_1(X^H)\rtimes \WH)$ and a pullback diagram  
\[\xymatrix{\overline{X_H}(x)\ar[r]^-{\tilde{\phi}}\ar[d]^{\pi_H{\tilde{\pi}_H}}&E(\pi_1(X^H)\rtimes \WH)\ar[d]^\pi\\X_H(x)/\WH\ar[r]^-\phi&B(\pi_1(X^H)\rtimes \WH)}\]
Then ${}_\phi B(\pi_1(X^H)\rtimes \WH)\odot B(\pi_1(X^H)\rtimes \WH)_\pi
\cong (X_H(x)/\WH)_{\pi_H{\tilde{\pi}_H}} \odot {}_ {\tilde{\phi}}E(\pi_1(X^H)\rtimes \WH)$.
Since $ {}_ {\tilde{\phi}}E(\pi_1(X^H)\rtimes \WH)\odot {}_{\pi\circ \tilde{i}^H}\ast\cong 
 {}_ {\tilde{\phi}\circ \pi\circ \tilde{i}^H}\ast$ and $ B(\pi_1(X^H)\rtimes \WH)_\pi\odot   {}_ {\tilde{\phi}\circ \pi\circ \tilde{i}^H}\ast\cong (E(\pi_1(X^H)\rtimes \WH),\pi)$ we can write $\widetilde{X^{H}}(x)\cup C(\rcov{X^{>H}}(x))$
as \[C_{H}X(x)\odot {}_\phi B(\pi_1(X^H)\rtimes \WH)\odot (E(\pi_1(X^H)\rtimes \WH),\pi).\]

As in the free case there is a corresponding decomposition of the map induced by $f$
and we have factoring of the trace.
\end{proof}

\section{Equivariant Nielsen numbers}\label{nn_sec}

Fadell and Wong \cite{FadellWong} and Wilczy\'nski \cite{Wilczynski}
have given very different proofs of the converse to the equivariant Lefschetz fixed point theorem. 
They used generalizations of the Nielsen number and we can compare their invariant with the equivariant Reidemeister trace 
using the results above.  We start with a consequence of Theorem A.

\begin{prop}\label{nrcompare}
$R_G(f)$ is zero if and only if $R(f^H)$ is zero for all subgroups $H$ of $G$. 
\end{prop}

Let $\Iso{X}$ be a choice of representatives for the isomorphism 
classes of objects in $\ecomcat{G}{X}$.  Without loss of generality we may assume that 
we have first chosen a representatives for each conjugacy class of subgroups of $G$ 
and that only these representatives appear among the objects of $\Iso{X}$.  Let $\Isog{X}{H}$ be the isomorphism 
classes of objects associated to maps $G/H\to X$.

\begin{lem}\label{rinject} For each subgroup $H$ of $G$ 
the forgetful map $\mathbb{Z}\Isog{X}{H}\to \mathbb{Z}\pi_0(X^H)$ defined by 
$x(H)\mapsto \sum_{gH\in \WH/H} x(gH)$ is injective.  
\end{lem}

\begin{proof}
The image of each $x(H)$ is nontrivial since all terms appear with coefficient 1.  

Suppose that $x(gH)$ and $y(g'H)$ are in the same component of $X^H$.  Then $x(H)$ and $g^{-1}g'y(H)$ are 
in the same component of $X^H$ and, in particular, $x(H)$ and $y(H)$ represent isomorphic objects in $\Isog{X}{H}$.  
\end{proof}

\begin{proof}[Proof of \autoref{nrcompare}]

The zeroth equivariant stable homotopy group of of a $G$-space $X_+$ 
is the free abelian group generated by $x(H)\circ \tr_G^\triangle (G/H)$ for  $x(H)\in \Iso{X}$ \cite[2.8.13.7]{tomDieck}. 
Combining coefficients if necessary, we can use Theorem A to express $R_G(f)$ as a sum \[\sum_{\gamma\in \Iso{\Lambda^fX}} a_\gamma (i_\gamma \circ \tr_G^\triangle (F_\gamma))\] 
where each element in $\Iso{\Lambda^fX}$ is associated to exactly one coefficient.  Then 
 $R_G(f)$ is zero if and only if each of the $a_\gamma$ are zero.

Additivity of the Reidemeister trace implies that $R(f^H)$ is zero for all subgroups $H$ of $G$ if and only if $R(f_H)$ is zero for all subgroups $H$ of $G$. 
Then Theorem A allows us to conclude that   \[R(f_{H})=\sum a_\gamma (i'_\gamma \circ \tr^\triangle(F'_\gamma))\]
where  we have exactly the same coefficients $a_\gamma$ as above and we sum over the elements of $\Isog{\Lambda^fX}{H}$. 

The invariant  $R(f_H)$ is the image of $R_{\WH}(f_H)$ under the forgetful map and so 
\autoref{rinject} implies $R(f_H)$ is zero if and only if 
each of the $a_\gamma$ are zero.
\end{proof}
 
Recall that two fixed points $x$ and $y$ in $X^H$ are  in the same \emph{fixed point class} if the images of the 
constant paths at $x$ and $y$ are in the same component of $\Lambda^{f^H}X^H$.  We say 
$x$ and $y$ are in the 
same $\WH$ \emph{fixed point class} if there is a $g\in \WH$ and so that the constant paths
at $x$ and $yg$ are in the same component of $\Lambda^{f^H}X^H$.  Note that there is a map from the 
fixed point classes to the $\WH$ fixed point classes.

If $H$ is subconjugate to $K$ the inclusion induces a map $\tau_{H\leq K}$ from 
the $\WK$ fixed point classes to the $\WH$ fixed point classes.
We let $\WH$-fpc denote the set of $\WH$ fixed point classes and
for $\alpha\in   \WH$-fpc let $i(f^H,\alpha)$ denote the nonequivariant
fixed point index of the fixed points in the class $\alpha$ with respect
to the map $f^H\colon X^H\to X^H$.  (At this point it is convenient to replace $f$ by a homotopic map 
that it taut and has isolated fixed points.  See \S\ref{symmetricmonoidal}.)

\begin{defn}\cite{Wong} 
  The \emph{equivariant Nielsen number} of $f$, $N_G(f)$, is the function from the 
  conjugacy classes of subgroups of 
  $G$ to the integers defined by 
  \[N_G(f)(H)=\sharp \left\lbrace \alpha\in \WH\text{-fpc}
  \left|
  \begin{array}{c}
  i\left(f^H,\alpha\right)\neq0 \text{ and }
  i\left(f^K,\delta\right)= 0\text{ for all } \\\delta\in \WK\text{-fpc}\text{ so that }
  \tau_{H\leq K}(\delta)=\alpha
  \end{array}
  \right\rbrace\right.\] 
\end{defn}

Note that the equivariant Nielsen number is \emph{not} the number of generators in 
the equivariant Reidemeister trace with nonzero coefficient.   The equivariant Nielsen number 
 is a  `non-redundant' count of the number of nonzero coefficients.  In particular, the 
coefficients of equivariant Reidemeister trace do not give a lower bound for the number of fixed points.

\begin{thm} The equivariant Nielsen number of a map is zero if and only if the equivariant Reidemeister trace is zero.
\end{thm}

\begin{proof}
For each $g\in \WH$ the diagram 
\[\xymatrix{
X^H\ar[r]^{(-)\cdot g}\ar[d]^f&X^H\ar[d]^f\\
X^H\ar[r]^{(-)\cdot g}&X^H
}\]
commutes and so for each fixed point $x$ of $f^H$, the indices of $x$ and $xg$ are the same \cite[IV.B.7]{brownbook}.  This implies that 
a $\WH$ fixed point class can only have index zero if all the associated (classical) fixed point classes have index zero.  

The equivariant Nielsen number is trivial if and only if the indices of all $\WH$-fixed point 
classes for all the induced maps $f^H\colon X^H\to X^H$ are zero.  
This holds if and only if the indices for all the fixed point classes 
of the maps $f^H$ are zero.  This is equivalent to  $R(f^H)$ is zero for all subgroups $H$ of $G$.  \autoref{nrcompare} completes the proof.
\end{proof}

It is unfortunate but necessary that the comparison in this section passes through the nonequivariant Reidemeister traces for isotropy subspaces.   
In particular,  both the equivariant Nielsen number and the equivariant Reidemeister trace can be used to give lower bounds for the 
number of fixed points, see \cite{Wong}  for the Nielsen number and \cite{gun}  for the Reidemeister trace, but this is lost with the classical Reidemeister trace.
This essential incompatibility is unsurprising since these two approaches are fundamentally very different.

\bibliographystyle{plain.bst}
\bibliography{trace}

\def\cprime{$'$}
\begin{thebibliography}{10}

\bibitem{Adams}
J.~F. Adams.
\newblock Prerequisites (on equivariant stable homotopy) for {C}arlsson's
  lecture.
\newblock In {\em Algebraic topology, Aarhus 1982 (Aarhus, 1982)}, volume 1051
  of {\em Lecture Notes in Math.}, pages 483--532. Springer, Berlin, 1984.

\bibitem{brownbook}
Robert~F. Brown.
\newblock {\em The {L}efschetz fixed point theorem}.
\newblock Scott, Foresman and Co., Glenview, Ill.-London, 1971.

\bibitem{crabb}
M.~C. Crabb.
\newblock The homotopy coincidence index.
\newblock {\em J. Fixed Point Theory Appl.}, 7(1):1--32, 2010.

\bibitem{daverman}
Robert~J. Daverman.
\newblock {\em Decompositions of manifolds}, volume 124 of {\em Pure and
  Applied Mathematics}.
\newblock Academic Press, Inc., Orlando, FL, 1986.

\bibitem{tomDieck}
Tammo~tom Dieck.
\newblock {\em Transformation groups}, volume~8 of {\em de Gruyter Studies in
  Mathematics}.
\newblock Walter de Gruyter \& Co., Berlin, 1987.

\bibitem{DP}
Albrecht Dold and Dieter Puppe.
\newblock Duality, trace, and transfer.
\newblock In {\em Proceedings of the International Conference on Geometric
  Topology (Warsaw, 1978)}, pages 81--102, Warsaw, 1980. PWN.

\bibitem{FadellWong}
Edward Fadell and Peter Wong.
\newblock On deforming {$G$}-maps to be fixed point free.
\newblock {\em Pacific J. Math.}, 132(2):277--281, 1988.

\bibitem{additivity}
Mortiz Groth, Kate Ponto, and Michael Shulman.
\newblock The additivity of traces in monoidal derivators.
\newblock arXiv:1212:3377, \emph{Journal of K-theory}.

\bibitem{KW2}
John~R. Klein and Bruce Williams.
\newblock Homotopical intersection theory. {II}. {E}quivariance.
\newblock {\em Math. Z.}, 264(4):849--880, 2010.
\newblock arXiv:0803.0017v2.

\bibitem{KW}
John~R. Klein and E.~Bruce Williams.
\newblock Homotopical intersection theory. {I}.
\newblock {\em Geom. Topol.}, 11:939--977, 2007.

\bibitem{Lait}
Erkki Laitinen and Wolfgang L{\"u}ck.
\newblock Equivariant {L}efschetz classes.
\newblock {\em Osaka J. Math.}, 26(3):491--525, 1989.

\bibitem{LMS}
L.~G. Lewis, Jr., J.~P. May, M.~Steinberger, and J.~E. McClure.
\newblock {\em Equivariant stable homotopy theory}, volume 1213 of {\em Lecture
  Notes in Mathematics}.
\newblock Springer-Verlag, Berlin, 1986.
\newblock With contributions by J. E. McClure.

\bibitem{LR}
Wolfgang L{\"u}ck and Jonathan Rosenberg.
\newblock The equivariant {L}efschetz fixed point theorem for proper cocompact
  {$G$}-manifolds.
\newblock In {\em High-dimensional manifold topology}, pages 322--361. World
  Sci. Publ., River Edge, NJ, 2003.

\bibitem{mayadd}
J.~P. May.
\newblock The additivity of traces in triangulated categories.
\newblock {\em Adv. Math.}, 163(1):34--73, 2001.

\bibitem{MS}
J.~P. May and J.~Sigurdsson.
\newblock {\em Parametrized homotopy theory}, volume 132 of {\em Mathematical
  Surveys and Monographs}.
\newblock American Mathematical Society, Providence, RI, 2006.

\bibitem{higher}
Kate Ponto.
\newblock Conicidence invariants and higher {R}eidemeister traces.
\newblock arXiv:1209.3710.

\bibitem{thesis}
Kate Ponto.
\newblock Fixed point theory and trace for bicategories.
\newblock {\em Ast\'erisque}, (333):xii+102, 2010.
\newblock arXiv:0807.1471.

\bibitem{relative}
Kate Ponto.
\newblock Relative fixed point theory.
\newblock {\em Algebr. Geom. Topol.}, 11(2):839--886, 2011.
\newblock arXiv:0906.0762.

\bibitem{PS_trace}
Kate Ponto and Michael Shulman.
\newblock Duality and trace in symemtric monoidal categories.
\newblock {\em Expositiones Mathematicae}, 2013.
\newblock 10.1016/j.exmath.2013.12.003. arXiv:1107.6032.

\bibitem{shadows}
Kate Ponto and Michael Shulman.
\newblock Shadows and traces in bicategories.
\newblock {\em J. Homotopy Relat. Struct.}, 8(2):151--200, 2013.
\newblock arXiv:0910.1306.

\bibitem{PS6}
Kate Ponto and Michael Shulman.
\newblock Linearity of fixed point invariants.
\newblock arXiv:1406.7861, 2014.

\bibitem{linearity}
Kate Ponto and Michael Shulman.
\newblock Linearity of traces in monoidal categories and bicategories.
\newblock arXiv:1406.7854, 2014.

\bibitem{PS_mult}
Kate Ponto and Michael Shulman.
\newblock The multiplicativity of fixed point invariants.
\newblock {\em Algebr. Geom. Topol.}, 14(3):1275--1306, 2014.
\newblock arXiv:1203.0950.

\bibitem{stallings}
John Stallings.
\newblock Centerless groups---an algebraic formulation of {G}ottlieb's theorem.
\newblock {\em Topology}, 4:129--134, 1965.

\bibitem{gun}
Gun Sunyeekhan.
\newblock {\em Equivariant intersection theory}.
\newblock ProQuest LLC, Ann Arbor, MI, 2010.
\newblock Thesis (Ph.D.)--University of Notre Dame.

\bibitem{ulrich}
Hanno Ulrich.
\newblock {\em Fixed point theory of parametrized equivariant maps}, volume
  1343 of {\em Lecture Notes in Mathematics}.
\newblock Springer-Verlag, Berlin, 1988.

\bibitem{Weber}
Julia Weber.
\newblock The universal functorial equivariant {L}efschetz invariant.
\newblock {\em $K$-Theory}, 36(1-2):169--207 (2006), 2005.

\bibitem{Weber2}
Julia Weber.
\newblock Equivariant {N}ielsen invariants for discrete groups.
\newblock {\em Pacific J. Math.}, 231(1):239--256, 2007.

\bibitem{whitehead}
J.~H.~C. Whitehead.
\newblock Simplicial spaces, nuclei and $m$-groups.
\newblock {\em Proc. London Math. Soc.}, s2-45(1):243--327, 1939.

\bibitem{Wilczynski}
Dariusz Wilczy{\'n}ski.
\newblock Fixed point free equivariant homotopy classes.
\newblock {\em Fund. Math.}, 123(1):47--60, 1984.

\bibitem{Wong}
Peter Wong.
\newblock Equivariant {N}ielsen fixed point theory for {$G$}-maps.
\newblock {\em Pacific J. Math.}, 150(1):179--200, 1991.

\bibitem{Wong2}
Peter Wong.
\newblock Equivariant {N}ielsen numbers.
\newblock {\em Pacific J. Math.}, 159(1):153--175, 1993.

\end{thebibliography}

\end{document}